\documentclass[11pt, leqno, a4paper, fleqno]{amsart}
\usepackage{amssymb,amsthm,amsmath,amscd,amsfonts,bbm,mathrsfs}
\usepackage[all]{xy}
\def\noproof{{\unskip\nobreak\hfill\penalty50\hskip2em\hbox{}%
     \nobreak\hfill$\Box$\parfillskip=0pt%
     \finalhyphendemerits=0\par}}
\def\enddemo{\ifmmode\eqno\Box\else\noproof\vskip0.8truecm\fi}

\newtheorem{theo}{Theorem}[section]
\newtheorem{theorem}[theo]{Theorem}
\newtheorem{definition}[theo]{Definition}
\newtheorem{lemma/def}[theo]{Lemma/Definition}

\newtheorem{prop}[theo]{Proposition}

\newtheorem{remarks}[theo]{Remarks}

\newtheorem{lemma}[theo]{Lemma}

\newcommand{\lra}{\longrightarrow}

\DeclareMathOperator{\Spec}{Spec}

\DeclareMathOperator{\Coind}{Coind}
\DeclareMathOperator{\Ind}{Ind}

\DeclareMathOperator{\res}{res}
\DeclareMathOperator{\cor}{cor}

\DeclareMathOperator{\sign}{sign}
\DeclareMathOperator{\solhu}{SH}

\DeclareMathOperator{\ord}{ord}

\DeclareMathOperator{\id}{id}

\DeclareMathOperator{\ev}{ev}
\DeclareMathOperator{\Hom}{Hom}
\DeclareMathOperator{\Aut}{Aut}

\DeclareMathOperator{\pr}{pr}

\DeclareMathOperator{\supp}{supp}

\DeclareMathOperator{\wcdot}{\, \cdot\, }

\DeclareMathOperator{\Ker}{Ker}
\DeclareMathOperator{\Log}{Log}
\DeclareMathOperator{\Image}{Im}

\DeclareMathOperator{\Maps}{Maps}
\DeclareMathOperator{\Real}{Re}

\DeclareMathOperator{\Cl}{Cl}

\DeclareMathOperator{\Gal}{Gal}
\DeclareMathOperator{\Norm}{N}

\DeclareMathOperator{\Sym}{Sym}

\DeclareMathOperator{\GL}{GL}
\DeclareMathOperator{\PGL}{PGL}

\DeclareMathOperator{\incl}{incl}

\DeclareMathOperator{\rat}{rat}

\newcommand{\llangle}{{\langle\!\langle}}
\newcommand{\rrangle}{{\rangle\!\rangle}}

\newcommand{\fa}{{\mathfrak a}}
\newcommand{\fb}{{\mathfrak b}}

\newcommand{\ff}{{\mathfrak f}}

\newcommand{\fm}{{\mathfrak m}}

\newcommand{\fp}{{\mathfrak p}}
\newcommand{\fq}{{\mathfrak q}}

\newcommand{\bC}{{\mathbb C}}

\newcommand{\bK}{{\mathbb K}}

\newcommand{\bN}{{\mathbb N}}
\newcommand{\bP}{{\mathbb P}}
\newcommand{\bQ}{{\mathbb Q}}
\newcommand{\bR}{{\mathbb R}}

\newcommand{\bZ}{{\mathbb Z}}
\newcommand{\barQ}{{\overline{\mathbb Q}}}

\newcommand{\bA}{{\mathbf A}}
\newcommand{\bI}{{\mathbf I}}

\newcommand{\barE}{{\overline{E}}}

\newcommand{\barcO}{{\overline{\cO}}}

\newcommand{\bfP}{{\bf P}}

\newcommand{\cA}{{\mathcal A}}

\newcommand{\cC}{{\mathcal C}}
\newcommand{\cD}{{\mathcal D}}

\newcommand{\cF}{{\mathcal F}}
\newcommand{\cG}{{\mathcal G}}
\newcommand{\cH}{{\mathcal H}}
\newcommand{\cI}{{\mathcal I}}
\newcommand{\cK}{{\mathcal K}}
\newcommand{\cL}{{\mathcal L}}

\newcommand{\cN}{{\mathcal N}}
\newcommand{\cO}{{\mathcal O}}
\newcommand{\cP}{{\mathcal P}}
\newcommand{\cQ}{{\mathcal Q}}
\newcommand{\cR}{{\mathcal R}}
\newcommand{\cS}{{\mathcal S}}

\newcommand{\wP}{\widetilde{P}}

\newcommand{\uv}{{\underline{v}}}
\newcommand{\ux}{{\underline{x}}}

\newcommand{\ep}{{\epsilon}}

\newcommand{\bu}{{\bullet}}

\newcommand{\noi}{\noindent}

\begin{document}

\title[Shintani cocycles and vanishing order of $L_p(\chi,s)$ at $s=0$]{Shintani cocycles and vanishing order of $p$-adic Hecke $L$-series at $s=0$}

\author{By Michael Spie{\ss}}
\date{March 27, 2012}
\address{Fakult\"{a}t f\"{u}r Mathematik,
Universit\"{a}t Bielefeld, D-33501 Bielefeld, Germany}
\email{mspiess@math.uni-bielefeld.de}
\begin{abstract}
Let $\chi$ be a Hecke character of finite order of a totally real number field $F$. By using Hill's Shintani cocyle we provide a cohomological construction of the $p$-adic $L$-series $L_p(\chi, s)$ associated to $\chi$. This is used to show that $L_p(\chi, s)$ has a trivial zero at $s=0$ of order at least equal to the number of places of $F$ above $p$ where the local component of $\chi$ is trivial. 
\end{abstract}
\maketitle
\tableofcontents

\section{Introduction} Let $F$ denote a totally real number field of degree $d>1$ over $\bQ$, let $p$ be a prime number and let $\chi$ be a totally odd Hecke character of finite order of $F$. Klingen and Siegel have shown that the values of the Hecke $L$-series $L(\chi,s)$ at integers $n\le 0$ lie in the algebraic closure $\barQ\subseteq \bC$ of $\bQ$. In \cite{shintani} Shintani gave another proof by constructing a nice fundamental domain (i.e.\ a finite disjoint union of {\it rational cones}; a so-called {\it Shintani decomposition}) for the canonical action of the positive global units $E_+$ of $F$ on $\bR_+^d$. 

Deligne and Ribet \cite{deligneribet} and independently Barsky and Cassou-Nogu{\`e}s \cite{barsky, cassou} have shown that there exists a $p$-adic analytic analogue $L_p(\chi,s)$ of the Hecke $L$-series $L(\chi, s)$ which is characterized by \footnote{This normalizations is not standard; usually this $p$-adic $L$-series is denote by $L_p(\chi\omega,s)$.}
\begin{equation*}
L_p(\chi, 1-n)\,\,\, = \,\,\, L_{S_p}(\chi\omega^{1-n}, 1-n)
\end{equation*}
for all integers $n\ge 1$. Here $\omega: \Gal(F(\mu_{2p})/F) \to (\bZ/2p\bZ)^* \to \bZ_p^*$ denotes the Teichm{\"u}ller character and $L_{S_p}(\chi, s)$ the $L$-series without the Euler factors at the places $F$ above $p$. Compared to Deligne-Ribet's work Barsky's and Cassou-Nogu{\`e}s' construction of $L_p(\chi,s)$ is more elementary and is based on Shintani's approach to the Theorem of Siegel-Klingen. 

Since $L(\chi, 0)\ne 0$ it follows in particular that $L_p(\chi, s)$ has a {\it trivial zero} at $s=0$ if and only if there exists a place $\fp$ above $p$ such that the local component $\chi_\fp$ of $\chi$ at $\fp$ is trivial. In \cite{gross} Gross conjectured that the order of vanishing  $\ord_{s=0} L_p(\chi, s)$ is equal to the number of places $\fp$ of $F$ above $p$ such that $\chi_\fp = 1$. 

In section 3 of our work \cite{ich} we developed a framework to deal with trivial zeros of higher order of $p$-adic $L$-functions. The latter are typically defined as so-called $\Gamma$-transforms of a $p$-adic measure on the Galois group of a certain infinite abelian extension $M/F$. In \cite{ich} we attach such a measure $\mu_{\kappa}$ to a cohomology class $\kappa\in H^{d-1}(F^*_+, \cD^b)$ where $F^*_+$ are the totally positive elements of $F$ and $\cD^b$ is a certain space of $p$-adic measures on the finite ideles $\prod_{v\nmid \infty}' F_v^*$ of $F$. We showed that the $\Gamma$-transform of $\mu_{\kappa}$ has a trivial zero of order at least $r$ (and also give a formula for its $r$-th derivative) if there exists $r$ places $\fp_1, \ldots, \fp_r$ of $F$ above $p$ such that $\kappa$ "extends" to a cohomology class whose values are measures on the larger adelic space $\prod_{i=1}^r F_{\fp_i} \times \prod_{v\nmid \infty, v \ne \fp_1, \ldots, \fp_r}' F_v^*$ (we will recall the set-up and results of \cite{ich} which are used in this paper in section \ref{section:ich} below).

In \cite{ich} we have applied this result to prove a conjecture of Hida regarding trivial zeros of the $p$-adic $L$-function $L_p(E, s)$ of a modular elliptic curve $E/F$. The aim of this paper is to apply it to $L_p(\chi,s)$ i.e.\ we give a proof of the following theorem.

\begin{theorem}
\label{theorem:vanishing}
Let $r$ be the number of places $\fp$ of $F$ above $p$ such that $\chi_\fp = 1$. Then, 
\begin{equation}
\label{vanorder}
\ord_{s=0} L_p(\chi, s) \ge r.
\end{equation}
\end{theorem}

We will work with Barsky's and Cassou-Nogu{\`e}s' construction of $L_p(\chi, s)$. However we need to "lift" the $p$-adic measure $\mu_{\chi}$ involved to a measure-valued cohomology class $\kappa_{\chi}$ in order to apply the method of \cite{ich}. This is achieved using a {\it Shintani cocyle}. It is a certain $(d-1)$-cocyle on $F^*_+$ with values in the module generated by all characteristic functions of rational cones in $\bR^d_+$ which yields a Shintani decomposition when taking the cap-product of it with the fundamental class in $H_{d-1}(E_+, \bZ)$ (for the precise definition see \ref{definition:shintani}). The notion of a Shintani cocyle has been introduced by Solomon \cite{solomon} who has given a definition in the case $d=2$. For arbitrary $d$, Hill \cite{hill} has given a construction. In section \ref{section:shintani} we recall it and -- by using a result of Colmez \cite{colmez} -- establish the relation to a Shintani decomposition (see Prop.\ \ref{prop:shincoc}). Then in section \ref{section:cassoutrick} and the beginning of section \ref{section:padiclseries} we carry out the construction of $\kappa_{\chi}$ which is followed by our proof of Thm.\ \ref{theorem:vanishing} (following Cassou-Nogu{\`e}s we choose a certain auxiliary prime $\fq\nmid p$ of $F$ to obtain $p$-integrality properties of twisted $L$-values). 

It should be mentioned that \eqref{vanorder} can be proved rather easily for the corresponding arithmetic $p$-adic $L$-function (see \cite{fedgross} and also \cite{burns}) and so Thm.\ \ref{theorem:vanishing} is a consequence of Iwasawa's deep {\it main conjecture} (as proven by Wiles \cite{wiles}). However we think that our approach is of independent interest. It is certainly more elementary. We also feel that the cohomological framework developed here might be useful to study other properties of $p$-adic $L$-series (we hope to return to the topic in the future). 

It should be mentioned as well that Dasgupta \cite{dasgupta2} independently (and earlier) gave a proof of \eqref{vanorder} if $r\le 3$ which is closely related to our approach. Moreover in joint work with Charollois \cite{daschar} he gives another proof of \eqref{vanorder} based in part on a different cohomological construction of $L_p(\chi, s)$ (involving Szech's {\it Eisenstein cocyle}).

While working on this paper I had helpful conversations with Pierre Colmez and Samit Dasgupta regarding "Cassou-Nogu{\`e}s' trick" so I thank them both.

\paragraph{\bf Notation} We introduce the following notation which will be used throughout the rest of this paper. 

We fix once and for all an embedding $\barQ\to \bC_p$. 

If $X$ and $Y$ are topological spaces then $C(X,Y)$ denotes the set of continuous maps $X\to Y$. If $R$ is a topological ring we let $C_c(X,R)$ denote the subset $C(X,R)$ of continuous maps with compact support. If we consider $Y$ (resp.\ $R$) with the discrete topology then we shall also write $C^0(X,Y)$ (resp.\ $C^0_c(X,R)$) instead of $C(X,Y)$ (resp.\ $C_c(X,R)$).

If $X$ is a locally compact Hausdorff space and $R = \bC_p$ we denote by $\|\wcdot \|_p$ the $p$-adic maximums norm on $C^0_c(X,\bC_p)$. It is given by
\begin{equation}
\label{maxnorm}
\| \phi\|_p \,\,\,=\,\,\, \max\{|\phi(x)|_p\mid \, x\in X\}\qquad \forall \, \phi\in \cC_c^0(X, \bC_p).
\end{equation}

For a group $G$ a subgroup $H$ there exists morphisms of $\delta$-functors
\[
\res: H^{\bu}(G, \wcdot) \lra H^{\bu}(H, \wcdot), \quad  \cor: H_{\bu}(H, \wcdot) \lra H_{\bu}(G, \wcdot).
\]
which in degree $0$ and for a $G$-module $M$ are the canonical inclusion $M^G\hookrightarrow M^H$ and projection $M_H \to M_G$ respectively. If $H$ has finite index in $G$ then there exists also morphisms
\[
\cor: H^{\bu}(H, \wcdot) \lra H^{\bu}(G, \wcdot), \quad  \res: H_{\bu}(G, \wcdot) \lra H_{\bu}(H, \wcdot).
\]
which in degree $0$ are given as follows. If $\{g_j\}_{j\in J}$ denotes systems of representatives of the cosets $G/H$ and $m\in M^H$ then $\cor(m) = \sum_{j\in J}\, g_j m$. For  $x = [m]\in M_G$ we have $\res(x) = [\sum_{j\in J}\, g_j^{-1} m]$.

Throughout the paper $F$ denotes a totally real number field of degree $d$ over $\bQ$ with ring of integers $\cO_F$. Let $E_F = \cO_F^*$ denote the group of global units. For a non-zero ideal $\fa\subseteq \cO_F$ we set $N(\fa) = \sharp(\cO_F/\fa)$. We denote by $\bfP_F$ the set of all places of $F$ and by $\bfP_F^{\infty}$ (resp.\ $S_{\infty}$) the subset of finite (resp.\ infinite) places. For a prime number $q$, we shall write $S_q$ for the set of places above $q$. We denote by $\sigma_1, \ldots, \sigma_d$ the different embeddings of $F$ into $\bR$.  Elements of $\bfP_F$ will be denoted by $v,w$ or also by $\fp, \fq$ if they are finite. If $\fp\in \bfP_F^{\infty}$, we denote the corresponding prime ideal of $\cO_F$ also by $\fp$. For $v\in \bfP_F$,  we denote by $F_v$ the completion of $F$ at $v$. If $v$ is finite then $\cO_v$ denotes the valuation ring of $F_v$ and $\ord_v$ the corresponding the normalized (additive) valuation on $F_v$ (so $\ord_v(\varpi) =1$ if $\varpi\in \cO_v$ is a local uniformizer at $v$). Also for $v\in \bfP_F$ we let $|\wcdot|_v$ be the associated normalize multiplicative valuation on $F_v$. Thus if $v\in S_{\infty}$ corresponds to the embedding $\sigma:F \to \bR$ then $|x|_v = |\sigma(x)|$ and if $v=\fq$ is finite then $|x|_{\fq} = N(\fq)^{-\ord_{\fq}(x)}$. For $v\in \bfP_F$ we put $U_v = \bR_+ = \{x\in \bR\mid x>0\}$ if $v$ is infinite and $U_v = \cO_v^*$ if $v$ is finite. 

Let $\bA= \bA_F$ be the adele ring of $F$ and $\bI= \bI_F$ the group of ideles. For a subset $S\subseteq \bfP_F$ we let $\bA\!^S$ (resp.\ $\bI^S$) denote the $S$-adeles (resp.\ $S$-ideles) and also define $\bA_S = \prod_{v\in S} F_v$ and $\bI_S= \prod_{v\in S} F_v^*$ . We also define $U^S =\prod_{v\not\in S} U_v$,  and $U_S = \prod_{v\in S} U_v$. If $S$ contains all archimedian places then the factor group $\bI^S/U^S$ is canonically isomorphic to the group $\cI^S$ of fractional $\cO_F$-ideals which are prime to all places in $S$. We sometimes view $F$ as a subring of $\bA_S$ and $\bA\!^S$ via the diagonal embedding. For a finite set of nonarchimedian places of $F$ we put $E^S = F^* \cap U^S$ and $E_S = F^* \cap U_S$ (intersection in $\bA_S$ resp.\ $\bA\!^S$).

For $T\subseteq \bfP_{\bQ} = \{ 2,3,5, \ldots, \infty\}$ and $S=\{v\in \bfP_F\mid\, v|_{\bQ}\in T\}$ we often write $\bA_T$, $\bA\!^T$, $\bI^T$ etc.\ for $\bA_S$, $\bA\!^S$, $\bI^S$ etc. We also write $U^q$, $U_q$, $U^{q, S}$, $U^{q,\infty}$ etc.\ for $U^{\{q\}}$, $U_{\{q\}}$, $U^{S_q\cup S}$, $U^{S_q\cup S_{\infty}}$ etc. and use a similar notation for adeles, ideles and fractional ideals. Thus for example for a finite subset $S$ of $\bfP_F^{\infty}$, $\bI^{S,\infty}$ denotes the set $S\cup S_{\infty}$-ideles and for a prime number $q$ we have $E_q = \{x\in F^*\mid \ord_v(x)= 0\,\forall\, v\,|\,q\}$. For $\ell\in \bfP_{\bQ} = \{ 2,3,5, \ldots, \infty\}$ we sometimes write $F_{\ell}$ rather than $\bA_{\ell} = F\otimes \bQ_{\ell} = \prod_{v\in S_{\ell}} F_v$. We shall denote by $F^*_+$, $E_{S,+}$, $E_+$ etc.\ the totally positive elements in $F$, $E_S$, $E$ etc.\ 

For a Hecke character $\chi: \bI/F^*\to \bC^*$ of finite order and $v\in \bfP_F$ we let $\chi_v$ be its $v$-component, i.e.\ $\chi_v : F_v^*\hookrightarrow \bI \stackrel{\chi}{\lra} \bC^*$. More generally if $S$ is a finite set of places of $F$ we denote $\chi_S: \bI_S\to \bC^*$ its restriction to $\bI_S\subseteq \bI$. If $S$ consists only of non-archimedian places the image of $\chi_S$ is contained in $\barQ^*\subseteq \bC_p$.

We denote by $\Norm = \Norm_{F/\bQ}: F^* \to \bQ^*$ the norm given by $\Norm(x) = \det(\ell_x)$ where $\ell_x:F \to F$ is the $\bQ$-linear map "multiplication by $x$". This extends to a map $(F\otimes_{\bQ} A)^* \to A^*$ for any $\bQ$-algebra $A$ which by abuse of notation will also be denoted by $\Norm$. 

Unless otherwise stated all rings are commutative with unit.

\section{$p$-adic $L$-series attached to cohomology classes}
\label{section:ich}

Let $R$ be a topological ring and let $S$ be a finite set of nonarchimedian places of $F$ (we mostly consider the case $S= S_p$). Assume that a decomposition $S = S_1 \overset{\cdot}{\cup} S_2$ into disjoint subsets is given. We introduce some spaces of $R$-valued functions  on adeles and ideles. Put $\cC(S_1,S_2, R) = C(\bA_{S_1}\times \bA_{S_2}^* \times \bI^{S,\infty}/ U^{S,\infty}, R)$, $\cC_c(S_1,S_2, R) = C_c(\bA_{S_1}\times \bA_{S_2}^* \times \bI^{S,\infty}/ U^{S,\infty}, R)$ and  $\cC^0_c(S_1, S_2, R) = C^0_c(\bA_{S_1}\times \bA_{S_2}^* \times \bI^{S,\infty}/U^{S,\infty}, R)$. We let $\bI^{\infty}$ act on $\cC_c(S_1, S_2, R)$ by $(a\cdot f)(x) = f(a^{-1}x)$ for $a\in \bI^{\infty}$, $f\in \cC_c(S_1, S_2, R)$ and $x\in \bA\!^{\infty}$. This induces an $F^*$- resp.\ $F_v^*$-action (for a finite place $v$) via the diagonal embedding $F^*\hookrightarrow \bI^{\infty}$ resp.\ the embedding $F_v^*\hookrightarrow \bI^{\infty}, x\mapsto (\ldots, 1, x, 1, \ldots)$. Note that 
\begin{equation*}
\label{unitint5}
\cC^0_c(S_1, S_2, R) \, \cong \, \bigotimes_{v\in S_1} C_c^0(F_v, R) \otimes \bigotimes_{v\in S_2} C_c^0(F_v^*, R) \otimes {\bigotimes_{v\not\in S}}' \, C_c^0(F_v,R)^{U_v}
\end{equation*}
where the restricted tensor product $\bigotimes'$ is taken with respect to the family of functions $\{1_{U_v}\}_v$. If $S_1 = \emptyset$ we often drop it from the notation, i.e.\ we write $\cC(S, R)$, $\cC_c(S, R)$ and $\cC^0_c(S, R)$ for $\cC(\emptyset, S, R) = C(\bI^{\infty}/U^{S,\infty}, R)$ etc.

Assume now that $S= S_p$ so $S_1\overset{\cdot}{\cup} S_2 = S_p$. Let $\cG_p= \Gal(M/F)$ be the Galois group of the maximal abelian extension $M/F$ which is unramified outside $p$ and $\infty$ and let $\rho: \bI/F^* \to \cG_p$ be the reciprocity map. There is a canonical homomorphism
\begin{equation}
\label{recipr2}
\partial: C(\cG_p, R) \lra  H_{d-1}(F^*_+, \cC_c(S_1, S_2,R))
\end{equation}
whose definition we recall from \cite{ich}. We denote by $\barE_+$ and $\overline{F^*_+}$ the closure of $E_+$ and $F^*_+$ with respect to the canonical embeddings
\begin{equation*}
\label{recipr2A}
E_+\hookrightarrow U_p,\quad F^*_+\hookrightarrow\bI^{\infty}/U^{p,\infty}
\end{equation*}
Note that $\overline{F^*_+} = F^*_+ \barE_+$. To begin with recall that the reciprocity map of class field theory $\rho: \bI/F^* \to \cG_p$ induces an isomorphism 
\begin{equation}
\label{recipr2B}
\bI^{\infty}/\overline{F_+^*}U^{p,\infty}    \,\,\, \cong \,\,\,   \bI/\overline{F^*}U^p          \,\,\, \cong \,\,\, \cG_p
\end{equation}
We can regard $\Gamma \colon = F^*_+/E_+$ as a discrete subgroup of $\bI^{\infty}/(\barE_+ \times U^{p,\infty})$ by using the embedding $\Gamma \hookrightarrow\bI^{\infty}/(\barE_+ \times U^{p,\infty})$. Next we construct an isomorphism
\begin{equation}
\label{recipr1}
H_0(\Gamma, \cC_c(S_p, R)^{E_+}) \lra C(\cG_p, R).
\end{equation}
Let $\pr: \bI^{\infty}/U^{p,\infty}\to \bI^{\infty}/(\barE_+ U^{p,\infty})$ denote the projection. Firstly, the map
\begin{equation}
\label{recipr2a}
C_c(\bI^{\infty}/\barE_+ U^{p,\infty}, R)\,\lra \, \cC_c(S_p, R)^{E_+},\,\,\, f\mapsto f\circ \pr
\end{equation}
is obviously an isomorphism. If $\cF$ denotes an (open and compact) fundamental domain of $\bI^{\infty}/\barE_+ U^{p,\infty}$ for the action of $\Gamma$ then 
\begin{equation*}
C_c(\bI^{\infty}/\barE_+ U^{p,\infty}, R) \, \cong \, \Ind^{\Gamma} \, C(\cF, R)
\, \cong \, \Ind^{\Gamma} \,C((\bI^{\infty}/\barE_+ \times U^{p,\infty})/\Gamma, R).
\end{equation*}
Hence there exists a canonical isomorphism
\begin{equation}
\label{recipr2b}
H_0(\Gamma, C_c(\bI^{\infty}/\barE_+ U^{p,\infty}, R))\lra C((\bI^{\infty}/\barE_+  U^{p,\infty})/\Gamma, R)
\end{equation}
which is given explicitly by $[f]\mapsto \sum_{[\zeta]\in \Gamma} \, f(\zeta x)$.  Since
\[
\left(\bI^{\infty}/\barE_+ \times U^{p,\infty}\right)/\Gamma \,\,\, \cong \,\,\, \bI^{\infty}/\overline{F_+^*}U^{p,\infty} \,\,\, \cong \,\,\,  \cG_p
\] 
(the second isomorphism is induced by the reciprocity map) the target of \eqref{recipr2b} can be identified with $C(\cG_p,R)$. Thus we have an isomorphism 
\begin{equation}
\label{recipr2c}
H_0(\Gamma, C_c(\bI^{\infty}/ \barE_+ \times U^{p,\infty}, R))\lra C(\cG_p, R)
\end{equation}
By combining \eqref{recipr2a} with the inverse of \eqref{recipr2c} we obtain the isomorphism \eqref{recipr1}.

Let $M$ be any $F^*_+$-module. Next we construct a homomorphism
\begin{equation}
\label{edge1}
H_0(\Gamma, H^0(E_+, M))\lra  H_{d-1}(F^*_+, M)
\end{equation}
Since $E_+ \cong \bZ^{d-1}$ we have $H_{d-1}(E_+,\bZ)\cong \bZ$. Choose a generator $\eta$ of $H_{d-1}(E_+,\bZ)$. Since the action of $\Gamma$ on $H_d(E_+,\bZ)$ is trivial, taking the cap product with $\eta$ yields an $\Gamma$-equivariant map $H^0(E_+, M) \to H_{d-1}(E_+, M)$ hence
\begin{equation}
\label{capeta}
H_0(\Gamma, H^0(E_+, M)) \lra  H_0(\Gamma, H_{d-1}(E_+, M)).
\end{equation}
We define \eqref{edge1} as the composite of \eqref{capeta} with the edge morphism 
\begin{equation*}
\label{edge2}
H_0(F^*_+/E_+, H_{d-1}(E_+, M))\lra H_{d-1}(F^*_+, M)
\end{equation*}
of the Hochschild-Serre spectral sequence. 

Finally we define \eqref{recipr2} to be the composite 
\begin{eqnarray*}
& C(\cG_p, R)\lra H_0(\Gamma, \cC_c(S_p, R)^{E_+}) \lra H_0(\Gamma, \cC_c(S_1 , S_2, R)^{E_+})\\
& \lra H_{d-1}(F^*_+, \cC_c(S_1 , S_2, R))
\end{eqnarray*}
where the first map is the inverse of \eqref{recipr1}, the second is induced by the natural inclusion $\cC_c(S, R) \hookrightarrow \cC_c(S_1 , S_2, R)$ and the third is \eqref{edge1} (for $M =  \cC_c(S_1 , S_2, R)$). If $R$ carries the discrete topology we write $\partial^0$ rather than $\partial$ for the map \eqref{recipr2}:
\begin{equation*}
\label{recipr2r}
\partial^0: C^0(\cG_p, R) \lra  H_{d-1}(F^*_+, \cC_c^0(S_1, S_2,R)).
\end{equation*}

\begin{remarks} 
\label{remarks:eta}
\rm (a) There is in fact a canonical choice for $\eta$ (since we have fixed a numbering $\sigma_1, \ldots ,\sigma_n$ of the embeddings $F\hookrightarrow \bR$). The norm $\Norm: F \to \bQ$ extends to a map $\Norm: F_{\infty}\to \bR$. We denote the kernel of $F_{\infty, +}^*\subseteq F_{\infty}^* \stackrel{\Norm}{\lra}\bR$ by $\cH$. The isomorphism $\Log: F_{\infty, +}^* \to \bR^d, \,y \mapsto (\log(\sigma_1(y)), \ldots, \log(\sigma_1(y))$ maps $\cH$ onto $\bR^d_0= \{z= (z_1,\ldots, z_d)\in \bR^d\mid \,  \sum_{i=1}^d z_i = 0\}$ and $E_+$ onto a complete lattice in $\bR^d_0$. The isomorphism $\cH/E_+\cong \bR^d_0/\Log(E_+)$ provides the $(d-1)$--dimensional compact manifold $\cH/E_+$ with an orientation. We chose $\eta\in H_{d-1}(E_+,\bZ)$ so that it corresponds  to the fundamental class under the canonical isomorphism $H_{d-1}(E_+,\bZ)\cong H_{d-1}(\cH/E_+,\bZ)$.

\noi (b) More generally for any discrete and cocompact subgroup of 
the $G$ of $\cH$ we obtain a canonical generator $\eta_G$ of $H_{d-1}(G, \bZ)$ in the same way as above. If $G'$ is a subgroup of finite index of $G$ then we have $\res(\eta_G) = \eta_{G'}$. 

\noi (c) The class $\eta_G$ can be described
explicitly in terms of generators $\ep_1, \ldots, \ep_{d-1}$ of $G$. Let $\mu=\pm 1$ be the sign of the determinant with rows $\Log(\ep_1), \ldots, \break \Log(\ep_{d-1}), v_0$  where $v_0 = (1, \ldots , 1)\in \bR^d$. Then $\eta_G$
 can be represented by the cycle $\mu\sum_{\tau\in S_{d-1}} \,\sign(\tau)\, [\ep_{\tau(1)}| \ldots| \ep_{\tau(d-1)}]$. 
 
 \noi (d) We shall also need the cohomological analogue of the map \eqref{edge1}. It is a homomorphism
\begin{equation}
\label{edgec1}
H^{d-1}(F^*_+, M)\lra H^0(\Gamma, H_0(E_+, M)).
\end{equation}
Taking the cap product with $\eta$ yields a $\Gamma$-equivariant map $H^{d-1}(E_+, \break M) \to H_0(E_+, M)$ hence
\begin{equation}
\label{capeta1}
H^0(\Gamma, H^{d-1}(E_+, M)) \lra  H^0(\Gamma, H_0(E_+, M))
\end{equation}
and \eqref{edgec1} is defined as the composite of the restriction
\begin{equation*}
\label{edgec2a}
H^{d-1}(F^*_+, M)\lra H^0(F^*_+/E_+, H^{d-1}(E_+, M))
\end{equation*}
followed by \eqref{capeta1}.
\end{remarks}

\paragraph{\bf Measure valued cohomology classes and $p$-adic $L$-functions} For a ring $R$ put
\[
\cD(S_1, S_2, R) = \Hom_R(\cC_c^0(S_1, S_2, R), R) = \Hom(\cC_c^0(S_1, S_2, \bZ), R)
\]
and let
\begin{equation}
\label{distr2}
\langle\,\,\, ,\,\,\, \rangle: \, \cD(S_1, S_2, R)\times \cC_c^0(S_1, S_2,R) \lra R,
\end{equation}
be the canonical (evaluation) pairing. We define a $\bI^{\infty}$-action on $\cD(S_1, S_2, R)$ by requiring that $\langle x\phi, x f\rangle = \langle \phi,  f\rangle$ for all $x\in \bI^{\infty}$, $f\in \cC^0_c(S_1, S_2, \break R)$ and $\phi\in \cD(S_1, S_2, R)$. The pairing \eqref{distr2} yields the bilinear map
\begin{eqnarray*}
\label{cap1}
&& \,\,\,\cap: H^{d-1}(F^*_+,\cD(S_1, S_2, R))\times H_{d-1}(F^*_+, \cC_c^0(S_1, S_2, R)) \lra \\
&& \hspace{2cm}\lra\,\, H_0(F^*_+, R) = R.\nonumber
\end{eqnarray*}

In the case $R = \bC_p$ an element $\lambda\in\cD(S_1, S_2, \bC_p)$ is called {\it bounded} if there exists a constant $C>0$ such that 
\begin{equation*}
\label{boundeddist}
|\lambda(\phi)|_p\,\,\, \le \,\,\,C\, \|\phi\|_p,\qquad \forall \phi\in \cC_c^0(S_1, S_2, \bC_p).
\end{equation*}
Here $\| \wcdot \|_p$ denotes the $p$-adic maximums norm  \eqref{maxnorm}. The $\bI^{\infty}$-submodule of bounded elements of $\cD(S_1, S_2, \bC_p)$ will be denoted by $\cD^b(S_1, S_2, \bC_p)$. Note that $\cD^b(S_1, S_2, \bC_p) = \Hom(\cC_c^0(S_1, S_2, \barcO_p), \barcO_p)\otimes_{\barcO_p} \bC_p$ where $\barcO_p$ denotes the valuation ring of $\bC_p$. Elements of $\cD(S_1, S_2, \bC_p)$ can be regarded as $\bC_p$-valued distributions and elements of $\cD^b(S_1, S_2, \bC_p)$ as $\bC_p$-valued measures on the locally compact space $\bA_{S_1}\times \bI^{S_1,\infty}/U^{p,\infty}$. Since a $\bC_p$-valued measure can be integrated against any continuous function with compact support the pairing \eqref{distr2} when restricted to $\cD^b(S_1, S_2, \bC_p)$ extends canonically to a pairing
\begin{equation*}
\label{meas2}
\langle\,\,\, ,\,\,\, \rangle: \,\cD^b(S_1, S_2, \bC_p)\times \cC_c(S_1, S_2,\bC_p) \lra \bC_p.
\end{equation*}
The latter gives rise to the bilinear map
\begin{equation*}
\label{cap3}
\cap: H^{d-1}(F^*_+,\cD^b(S_1, S_2, \bC_p))\times H_{d-1}(F^*_+, \cC_c(S_1, S_2, \bC_p)) \lra \bC_p.
\end{equation*}
For $\kappa\in H^{d-1}(F^*_+,\cD^b(S_1, S_2, \bC_p))$ we define the $\bC_p$-valued distribution $\mu_{\kappa}$ on $\cG_p$ by
\begin{equation}
\label{cap2a}
\int_{\cG_p} f(\gamma) \,\mu_{\kappa}(d\gamma) \,\, = \,\, \iota_*(\kappa)\cap \partial^0(f) \hspace{2cm} \forall\,\, f\in C^0(\cG_p, R)
\end{equation}
(here $\iota: \cD^b(S_1, S_2, \bC_p) \hookrightarrow \cD(S_1, S_2, \bC_p)$ denotes the inclusion). 

\begin{lemma} 
\label{lemma:distmeas} 
(a) $\mu_{\kappa}$ is a $p$-adic measure on $\cG_p$.

\noi (b) $\int_{\cG_p} f(\gamma) \,\mu_{\kappa}(d\gamma) \, =\, \kappa \cap \partial(f)$ for all $f\in C(\cG_p, \bC_p)$.
\end{lemma}

{\em Proof.} It suffices to show that there exists a constant $C>0$ with
\begin{equation}
\label{meas3}
\hspace{0.5cm}|\kappa \cap \partial(f)|_p \,\,\le \,\, C\, \|f\|_p  \hspace{2cm}\forall \, f\in C(\cG_p, \bC_p).
\end{equation}
Choose $\lambda\in \cD^b(S_p, \bC_p)$ representing the class of the image of $\kappa$ under the map \eqref{edgec1} i.e.\ 
\(
H^{d-1}(F^*_+,\cD^b(S_p, \bC_p))\to H^0(\Gamma, H_0(E_+, \cD^b(S_p, \bC_p)))
\).
Also let $\cF$ denote an open and compact subset of $\bI^{\infty}$ which is $U^{\infty}$-stable and such that the image of $\cF$ under $\pr: \bI^{\infty}\to \bI^{\infty}/U^{\infty}$ is a fundamental domain for the action of $\Gamma = F^*_+/E_+$. Then \eqref{meas3} follows immediately from
\[
\hspace{1cm}\kappa \cap \partial(f) \,\, =\,\, \lambda(1_{\cF} \cdot f\circ \rho)\hspace{2cm}\forall \, f\in C(\cG_p, \bC_p)
\]
where $\rho: \bI^{\infty}/U^{p, \infty} \to \cG_p$ denotes the reciprocity map. 
\enddemo

Recall that the $\Gamma$-transform of a $\bC_p$-valued $p$-adic measure $\mu$ on $\cG_p$ is given by 
\begin{equation*}
\label{padiclfunc}
L_p(\mu, s) \, =\, \int_{\cG_p}\, \langle \gamma\rangle^s \mu(d\gamma)
\end{equation*}
for $s\in \bZ_p$. Here $\langle \gamma\rangle^s \colon = \exp_p(s \log_p(\cN(\gamma)))$ where $\cN: \cG_p \to \bZ_p^*$  is the cyclotomic character (it is characterized by $\gamma \zeta = \zeta^{\cN(\gamma)}$ for all $p$-power roots of unity $\zeta$). We have (\cite{ich}, Thm.\ 3.13)

\begin{theorem}
\label{theorem:abstrezc2} 
Let $S_1, S_2$ be disjoint subsets of $S_p$ with $S_1\cup S_2 = S_p$, let $\kappa\in H^{d-1}(F^*_+, \cD^b(S_1,  S_2, \bC_p))$ and let $\mu = \mu_{\kappa}$ be the associated $p$-adic measure on $\cG_p$. Then, 
\[
\ord_{s=0} L_p(\mu, s) \ge \sharp(S_1).
\]
\end{theorem}

We point out that this is consequence of $\partial((\log_p \circ \cN)^k)=0$ for $k=0, \ldots r-1$ (where $r= \sharp(S_1)$), a fact which is proved in \cite{ich}. Indeed, by \ref{lemma:distmeas} (b) we get for the $k$-th derivative of $L_p(\kappa, s)$ at $s=0$ 
\begin{equation*}
L^{(k)}_p(\mu, 0) \,=\, \int_{\cG_p}\, (\log_p\circ \cN)^k \mu(d\gamma)\,=\, \kappa \cap \partial((\log_p\circ \cN)^k) 
\end{equation*}
hence $\ord_{s=0} L_p(\mu, s) \ge r$.

\paragraph{\bf Variant} Assume again that a decomposition $S = S_1 \overset{\cdot}{\cup} S_2$ of an arbitrary finite subset $S$ of $\bfP_F^{\infty}$ into disjoint subsets is given. For a finite set $T\subseteq \bfP_F^{\infty}$ disjoint to $S$ we define $\cC(S_1,S_2, R)^T$ (resp.\ $\cC(S_1,S_2, R)_T$) by omitting all places of $T$ (resp.\ by omitting all places of $F$ not lying in $S\cup T$) from the definition of $\cC(S_1,S_2, R)$. Thus
\begin{eqnarray*}
\cC(S_1,S_2, R)^T & = & C(\bA_{S_1}\times \bA_{S_2}^* \times \bI^{S\cup T, \infty}/ U^{S\cup T, \infty}, R)\\
\cC(S_1,S_2, R)_T & = & C(\bA_{S_1}\times \bA_{S_2}^* \times \bI_T/ U_T, R).
\end{eqnarray*}
Similarly one defines $\cC^0_c(S_1, S_2, R)^T$, $\cC^0_c(S_1, S_2, R)_T$, $\cD(S_1, S_2, R)^T$,\newline $\cD(S_1, S_2, R)_T$, $\cD^b(S_1, S_2, \bC_p)^T$,  etc.\ If $T = S_q$ for a prime number $q$ then we shall also write $\cC(S_1,S_2, R)^q$, $\cC^0_c(S_1, S_2, R)^q$ etc.\ for $\cC(S_1,S_2, R)^T$, $\cC^0_c(S_1, S_2, \break R)^T$ etc.\

The $R$-module $\cC^0_c(S_1, S_2, R)^T$ (resp.\ $\cC^0_c(S_1, S_2, R)_T$) carries a canonical $\bI^{T, \infty}$-action (resp.\ $\bI_T$-action). It is easy to see that 
\begin{eqnarray*}
\label{shapiro}
&&\cC(S_1,S_2, R) \,\,\, \cong \,\,\, \Ind^{\bI^{\infty}}_{\bI^{T, \infty}} \cC(S_1,S_2, R)^T\\
&&\cC(S_1,S_2, R) \,\,\, \cong \,\,\, \Ind^{\bI^{\infty}}_{\bI_T} \cC(S_1,S_2, R)_T.\nonumber
\end{eqnarray*}
Thus by weak approximation we have $\cC(S_1,S_2, R)\cong \Ind^{F_+^*}_{E_{T,+}} \cC(S_1,S_2, R)^T$ as $F_+^*$-modules. Moreover if $\bI^{S\cup T, \infty} = U^{S\cup T, \infty}F_+^*$ then
\begin{equation*}
\cC(S_1,S_2, R) \,\,\, \cong \,\,\, \Ind^{F_+^*}_{E^{S\cup T}_+} \cC(S_1,S_2, R)_T.
\end{equation*}
Similar statements hold for $\cC_c(S_1, S_2, q, R)$, $\cD(S_1, S_2, R)$ and $\cD^b(S_1, S_2, \bC_p)$ where in the latter cases $\Ind$ has to be replaced by $\Coind$.  

For example by Shapiro's Lemma there are canonical isomorphisms
\begin{eqnarray}
&& H_{d-1}(E^{S\cup T}_+,\cC_c(S_1, S_2, R)_T) \,\,\, \cong \,\,\, H_{d-1}(F^*_+, \cC_c(S_1, S_2, R))
\label{shapiro1}\\
&& H^{d-1}(E^{S\cup T}_+,\cD^b(S_1, S_2, \bC_p))_T \,\,\, \cong \,\,\, H^{d-1}(F^*_+,\cD^b(S_1, S_2, \bC_p))\nonumber
\end{eqnarray}
if $\bI^{S\cup T, \infty} = U^{S\cup T, \infty}F_+^*$.

Hence the map \eqref{recipr2}, Thm.\ \ref{theorem:abstrezc2} etc.\ can be formulated in terms of $E_{T,+}$-(co-)homology with coefficients in $\cC^0_c(S_1, S_2, R)_T$, $\cD(S_1, S_2, \bC_p)_T$ etc.\ The isomorphisms \eqref{shapiro} are induced by inclusions (i.e.\ adjunction maps)
\begin{eqnarray}
\label{shapiro2}
\cC(S_1,S_2, R)^T & \lra & \cC(S_1,S_2, R)\\
\label{shapiro3}
\cC(S_1,S_2, R)_T & \lra & \cC(S_1,S_2, R) 
\end{eqnarray}
given by $\phi \mapsto \phi\otimes 1_{U_T}$ (resp.\ $\phi \mapsto \phi\otimes 1_{U^{S\cup T, \infty}}$). Note that \eqref{shapiro2} is $E_T$-equivariant and  \eqref{shapiro3} is $E^{S\cup T}$-equivariant. 

If $q$ is a prime number such that $S_q$ is disjoint from $S\cup T$ then \eqref{shapiro3} factors in the form $\cC(S_1,S_2, R)_T \stackrel{\eqref{shapiro4}}{\lra}  \cC(S_1,S_2, R)^q \stackrel{\eqref{shapiro2}}{\lra} \cC(S_1,S_2, R)$ where 
\begin{equation}
\label{shapiro4}
\cC(S_1,S_2, R)_T\,\,\lra \,\, \cC(S_1,S_2, R)^q, \,\, \phi \mapsto \,\, \phi\otimes 1_{U^{S\cup T, q , \infty}}.
\end{equation}

\section{Shintani cocycles}
\label{section:shintani}

\paragraph{\bf Definition of a Shintani cocyle} For linearly independent $x_1, \ldots x_m \in F_{\infty}$ the (open) cone $C(x_1,  \ldots, x_m)$ generated by $(x_1, \ldots x_m)$ is defined by 
\[
C(x_1,  \ldots x_m)\,\, =\,\, \left\{\sum_{i=1}^m \, \lambda_i x_i\mid \, \lambda_i\in \bR_+  \,\,\forall \, i=1, \ldots ,m\right\}.
\]
It is called positive if $x_1, \ldots x_m\in F_{\infty, +}$ and rational if it is generated by $x_1, \ldots x_m\in F^*$. A rational cone $C= C(x_1,  \ldots x_m)$ generated by $x_1, \ldots x_m\in F_+$ is called a {\it Shintani cone}. Note that this is equivalent to the condition $C\subseteq F_{\infty, +}$. A subset $D$ of $F_{\infty, +}$ that can be written as a finite disjoint union of Shintani cones is called a {\it Shintani set}. Let $\cK$ (resp.\ $\cK_{\rat}$) denote the $\bZ$-span generated by all characteristic functions $1_C$ of positive (resp.\ Shintani) cones. The group $F_{\infty, +}^*$ acts on $\cK$ via $(x \cdot f)(y) = f(x^{-1} y)$ for $x\in F_{\infty, +}^*$, $f\in \cK$ and $y\in F_{\infty}$ and $\cK_{\rat}$ is a $F_+^*$-stable subgroup. Since the intersection of two Shintani sets is again a Shintani set the product of two functions in $\cK_{\rat}$ lies again in $\cK_{\rat}$.

Let $G$ be a discrete subgroup of $\cH = \{x\in F_{\infty, +}^*\mid \Norm(x) =1 \}$. It intersects a given positive cone $C = C(x_1,  \ldots x_m)$ in only finitely many points. In fact $\Log$  maps $G$ to a lattice in $\bR^d_0$ and $C\cap \cH$ to a bounded open subset of $\bR^d_0$ so the intersection $G\cap C \approx \Log(G \cap C)$ is finite. Hence for $f\in \cK$ and $y\in F_{\infty,+}$ almost all terms in the sum $\sum_{x\in G}\, (x f)(y)$ are $=0$. Thus the map $f\mapsto \sum_{x\in G}\, x f$ induces a homomorphism
\begin{equation}
\label{shinpair}
H_0(G, \cK) \lra \Maps(F_{\infty,+}, \bZ)^G.
\end{equation}
For a subgroup $G$ of $\cH$ which is discrete and cocompact (i.e.\ $\Log(G)$ is a complete lattice in $\bR^d_0$) we let 
\begin{equation}
\label{shinmap}
\psi_G: H^{d-1}(\cH, \cK)\lra  \Maps(F_{\infty,+}, \bZ)
\end{equation}
be the composition of 
\begin{equation*}
\label{edgec2b}
H^{d-1}(\cH, \cK)\stackrel{\res}{\lra} H^{d-1}(G, \cK)\stackrel{\cap \eta_{G}}{\lra} H_0(G, \cK) 
\end{equation*}
with the map \eqref{shinpair}.

\begin{lemma} 
\label{lemma:rescor} 
Let $G, G_1, G_2$ be discrete cocompact subgroups of $\cH$. 

\noi (a) If $G_1\subseteq G_2$ then $\psi_{G_1} = \psi_{G_2}$. 

\noi (b) Put $\sqrt[\infty]{G} \colon = \{x\in F_{\infty, +}^*\mid x^n \in G \,\, \mbox{for some $n\in \bN$}\,\}$. Then the image of \eqref{shinmap} is $\sqrt[\infty]{G}$-invariant.
\end{lemma}

{\em Proof.} (a) This follows from the commutativity of the diagram
\begin{equation*}
\begin{CD}
H^{d-1}(G_2, \cK)@>\cap \eta_{G_2} >> H_0(G_2, \cK) @>  \eqref{shinpair} >> \Maps(F_{\infty,+}, \bZ)^{G_2}\\
@VV\res V @VV \res V@VV\incl V  \\
H^{d-1}(G_1, \cK)@> \cap \eta_{G_1} >> H_0(G_1, \cK) @>  \eqref{shinpair} >> \Maps(F_{\infty,+}, \bZ)^{G_1}
\end{CD}
\end{equation*}
which is due to $\res\eta_{G_2} = \eta_{G_1}$. 

\noi (b) is a consequence of (a) since $\sqrt[\infty]{G}$ is the union of groups $G'$ with $G' \supseteq G$ and $[G':G]<\infty$ and since $\Image(\psi_{G'})$ is $G'$-invariant.
\enddemo

In the case $G=E_+$ we restrict the map $f\mapsto \sum_{\ep\in E_+}\, \ep f$ to $f\in \cK_{\rat}$ and obtain a homomorphism 
\begin{equation}
\label{shinpair2}
H_0(E_+, \cK_{\rat}) \lra \Maps(F_{\infty,+}, \bZ)^{E_+}.
\end{equation}
Recall that a Shintani set $\cA$ is called a {\it Shintani decomposition} if $F_{\infty,+}$ can be decomposed as the disjoint union of $E_+$-translates of $\cA$
\begin{equation*}
\label{shindecomp}
F_{\infty,+} \,\,\, =\,\,\, \overset{\cdot}{\bigcup}_{\epsilon \in E'} \ep \cA.
\end{equation*} 
Shintani \cite{shintani} has shown that they exists.

\begin{lemma}
\label{lemma:shincocycel1}
The map \eqref{shinpair2} is injective.
\end{lemma}

{\em Proof.} Let $f\in \cK$ and let $\cA$ be a Shintani decomposition for $E_+$. By (\cite{dasgupta1}, Lemma 3,14) there exists only finitely many $\ep\in E_+$ with $1_{\ep\cA} \cdot f \ne 0$. Modulo the augmentation ideal of the group ring $\bZ[E_+]$ we get
\begin{equation*} 
f \,\, = \,\, \sum_{\ep\in E_+}\, 1_{\ep\cA} \cdot f \,\, =\,\,\sum_{\ep\in E_+}\, \ep (1_{\cA} \cdot \ep^{-1} f) \,\, \equiv \,\, \sum_{\ep\in E_+}\, 1_{\cA}\cdot (\ep^{-1} f).
\end{equation*}
Hence any element of $H_0(E_+, \cK)$ has a representative $f\in \cK$ with $\supp(f) = \{x\in F_{\infty, +}\mid f(x) \ne 0\}\subseteq \cA$. For such $f$ and $y\in \cA$ we have $\sum_{\ep\in E_+}\, (\ep f)(y) = f(y)$ so in particular $\sum_{\ep\in E_+}\, \ep f=0$ implies $f=0$. \enddemo

Similarly to \eqref{shinmap} we define  
\begin{equation}
\label{shinmap2}
H^{d-1}(F_+^*, \cK_{\rat})\lra  \Maps(F_{\infty,+}, \bZ)
\end{equation}
to be the composition of 
\begin{equation}
\label{edgec2}
H^{d-1}(F_+^*, \cK_{\rat})\stackrel{\res}{\lra} H^{d-1}(E_+, \cK_{\rat})\stackrel{\cap \eta_{E_+}}{\lra} H_0(E_+, \cK_{\rat}) 
\end{equation}
with the map \eqref{shinpair2}.

\begin{definition}
\label{definition:shintani}
A cohomology class $\alpha\in  H^{d-1}(F_+^*, \cK_{\rat})$ is called Shintani cocycle if the following equivalent conditions hold:

\noi (i) $\alpha$ is mapped under \eqref{edgec2} to the class $[1_{\cA}]\in H_0(E_+, \cK_{\rat})$ of the characteristic function $1_{\cA}$ of a Shintani decomposition $\cA$. 

\noi (ii) $\alpha$ is mapped under \eqref{shinmap2} to the constant function $\equiv 1$ on $F_{\infty,+}$.
\end{definition}
 
That (i) and (ii) are indeed equivalent follows from Lemma \ref{lemma:shincocycel1} since $[1_{\cA}]$ is mapped to the constant function $\equiv 1$ under \eqref{shinpair2}.

Solomon \cite{solomon} has defined a Shintani cocyle in the case $d=2$. In \cite{hill} Hill has constructed a cohomology class in $H^{d-1}(F^*, \cL)$ for arbitrary $d$ where $\cL$ is the quotient of the $\bZ$-span of the characteristic functions of all rational cones modulo the constant functions $F_{\infty}-\{0\}\to \bZ$. We recall his construction and show that it is indeed a Shintani cocyle in the above sense when restricted to the subgroup $F_+^*$ of $F^*$.

\paragraph{\bf Hill's construction} Let $k$ be a field and $V$ a $k$-vector space of dimension $d$. Let $K/k$ be a field extension and let $t_1, \ldots, t_d\in K$ be elements which are algebraically independent over $k$. For a $k$-basis $\uv = (v_1, \ldots, v_d)$ of $V$, a $k$-algebra $A$ and $a\in A$ we define the element $b(\uv, a)\in V_A= A\otimes_k V$ by 
\begin{equation*}
\label{hillvector} 
b(\uv, a) \,\, =\,\, \sum_{j=1}^d \, a^{i-1}\, v_i.
\end{equation*}
Recall that a set of vectors of $V$ is said to be {\it in general position} if every subset with no more than $d$ vectors is linearly independent. We have (see also \cite{hill}, Lemma 1)

\begin{lemma}
\label{lemma:shin1}
Let $K/k$ be a field extension and let $t_1, \ldots, t_d\in K$ be elements which are algebraically independent over $k$. Given a collection of $d$ bases $\uv^{(1)} = (v_1^{(1)}, \ldots , v_d^{(1)})$, \ldots, $\uv^{(d)}= (v_1^{(d)}, \ldots , v_d^{(d)})$ of $V$ and $v\in V$, $v\ne 0$ the set of vectors in $V_K$
 \[
(v, b(\uv^{(1)}, t_1), \ldots, b(\uv^{(d)}, t_d))
\]
is in general position.
\end{lemma}

{\em Proof.} Let $\omega: V^d\to k$ be a determinant form on $V$. To prove that $(b(\uv^{(1)}, t_1), \ldots, b(\uv^{(d)}, t_d))$ is linearly independent it suffices to show
\begin{equation*}
\label{det} 
f(t_1, \ldots, t_d) \colon = \omega(b(\uv^{(1)}, t_1), \ldots, b(\uv^{(d)}, t_d)) \,\, \ne \,\, 0.
\end{equation*}
Let $X$ be the closed subvariety $f=0$ of $\bA\!^n_k = \Spec k[t_1, \ldots, t_d]$. We have to show $X\ne \bA\!^n_k$. For that we may assume that $k$ is infinite. Now for any $d$ elements $a_1, \ldots, a_d\in k$ with $a_i \ne a_j$ for $i\ne j$ the $d$-tuple $(b(\uv, a_1), \ldots, b(\uv, a_d))$ is a basis of $V$. Hence by Steinitz' lemma we can choose $a_1, \ldots, a_d\in k$ successively so that $(b(\uv^{(1)}, a_1), \ldots, b(\uv^{(d)}, a_d))$ is a basis, i.e.\ $f(a_1, \ldots, a_d) \ne 0$. Similarly, one shows that $(v, b(\uv^{(1)}, t_1), \ldots, b(\uv^{(d-1)}, t_{d-1}))$ is linearly independent.\enddemo

Assume now that $k$ is an ordered field (i.e.\ $k=\bQ$ or $k= \bR$), $V$ a $d$-dimensional $k$-vector space and $\omega: V^d\to k$ a determinant form on $V$. For an extension of ordered fields $K/k$ and linearly independent vectors $v_1, \ldots, v_m$ of $V_K$ we denote as before the open cone $\{\sum_{i=1}^m \, \lambda_i v_i\mid \, \lambda_i\in K, \,\lambda_i >0 \,\forall \, i=1, \ldots ,m\}$ by $C(v_1,  \ldots v_m)$. We also use the notation $C[v_1,  \ldots v_m]$ for the  closed cone $\{\sum_{i=1}^m \, \lambda_i v_i\mid \, \lambda_i\in K, \, \lambda_i \ge 0 \,\,\forall \, i=1, \ldots ,m\}-\{0\}$. A cone is called $k$-rational if the generators $v_1, \ldots, v_m$ can be chosen in $V$. If $d=m$ we define a function $c(v_1,  \ldots v_d): V_K-\{0\}\to \bQ$ by
\begin{equation*}
\label{charcone} 
c(v_1, \ldots, v_d)(v) \,\, =\,\, \sign(\omega(v_1, \ldots, v_d))\, 1_{C(v_1,  \ldots v_d)}(v). 
\end{equation*}

We choose a field extension $K/k$ with a discrete rank $d+1$ valuation $w$ which is trivial on $k$ and so that its residue field is $k$. Recall that $w$ is a surjective map $w: K \to \bZ^{d+1}\cup\{+\infty\}$ such that for all $x,y\in K$ we have (i) $w(x) = +\infty\Leftrightarrow x=0$, (ii) $w(xy) = w(x) + w(y)$ and (iii) $w(x+y) \ge \min(w(x), w(y))$. Here $\bZ^{d+1}$ carries the lexicographical order $(m_0, \ldots, m_d)< (n_0, \ldots, n_d)\Leftrightarrow m_0 = n_0, \ldots, m_{i-1} = n_{i-1}, m_i< n_i$ for some $i\in \{0, \ldots, d\}$. We denote the corresponding valuation ring by $\cO_w$ and the valuation ideal by $\fm_w$ so that $\cO_w/\fm_w =k$. For $x\in \cO_w$ we denote by $\bar{x}\in k$ the residue class $x \!\!\!\mod \!\fm_w$. We fix $t_0, t_1, \ldots, t_d\in K$ so that $w(t_0) = (1,0,\ldots, 0), w(t_1) = (0,1,\ldots, 0), \ldots, w(t_d) = (0,\ldots, 0, 1)$ is the standard basis of $\bZ^{d+1}$. We introduce an ordering on $K$ by 
\begin{equation}
\label{order}
x \, > \, 0 \quad \Leftrightarrow \quad \bar{u} \, >\, 0 
\end{equation}
where $x =u t_0^{m_0} \cdots t_d^{m_d}\in K^*$ with $u\in \cO_w^*$ and $w(x) = (m_0, \ldots, m_d)$. Note that $t_0, \ldots, t_d$ are algebraically independent over $k$ and that we have $0 < t_d$ and $t_i < t_{i-1}^m$ for all $i\in \{1, \ldots, d\}$ and $m\in \bN$. 

We fix a $k$-basis $\uv = (v_1, \ldots, v_d)$ of $V$ and a determinant form $\omega: V^d\to k$ on $V$. For $g_1, \ldots, g_d\in \GL(V) \colon = \Aut(V)$ and $v\in V-\{0\}$ we define 
\begin{equation*}
\label{hillcocycle} 
z(g_1, \ldots, g_d)(v) \,\, =\,\, c(g_1(b(\uv, t_1)), \ldots, g_d(b(\uv, t_d)))(v).
\end{equation*}
Note that $g b(\uv, a) = b(g\uv, a)$ for all $g\in \GL(V)$ and $a\in K$. Hence by Lemma \ref{lemma:shin1} the vectors $v, g_1 b(\uv, t_1)), \ldots, g_d b(\uv, t_d))$ are in general position. By (\cite{hill}, Prop.\ 3 and Thm.\ 1) the following holds

\begin{itemize}
\item[(SC1)] For all $g_0, g_1, \ldots, g_d\in \GL(V)$ and $v\in V-\{0\}$ we have
\begin{eqnarray*}
\label{shinc1} 
&& \sum_{i=0}^d \, (-1)^i \,z(g_0, \ldots, \widehat{g_i},\ldots, g_d)(v) \\
&& \hspace{1cm}=\,\, (-1)^d \, c(b(g_0\uv, t_0), \ldots, b(g_{d-1}\uv, t_{d-1}))(-b(g_d\uv, t_d)).\nonumber
\end{eqnarray*}
\item[(SC2)] For $g, g_1, \ldots, g_d\in \GL(V)$ and $v\in V-\{0\}$ we have
\begin{equation*}
\label{shinc2} 
z(g g_1, \ldots, g g_d)(v) \,\, =\,\, \sign(\det(g)) \,z(g_1, \ldots, g_d)(g^{-1}v) \end{equation*}
\item[(SC3)] For $g_1, \ldots, g_d\in \GL(V)$ there are finitely many disjoint $k$-ratio\-nal open cones $C_1, \ldots, C_m$ such that $z(g_1, \ldots, g_d) = \sum_{i=1}^m \, 1_{C_i}$.
\end{itemize}
We also need the following 

\begin{lemma}
\label{lemma:shin2}
Let $g_1, \ldots, g_d\in \GL(V)$ and assume that $b_1= g_1(v_1), \ldots, b_d = g_d(v_1)$ is a basis of $V$. Then 
\[
\sign(\omega(b_1, \ldots, b_d)) \,= \,\sign(\omega(g_1(b(\uv, t_1)), \ldots, g_d(b(\uv, t_d)))
\]
and 
\[
C(b_1, \ldots,  b_d) \,\subseteq\, C(g_1(b(\uv, t_1)), \ldots, g_d(b(\uv, t_d))\cap V \,\subseteq\, C[b_1, \ldots,  b_d].
\]
\end{lemma}

{\em Proof.} The first assertion follows immediately from the definition of the ordering on $K$. Let $v\in V-\{0\}$ and let $\lambda_1, \ldots, \lambda_d\in K$ be the coordinates of $v$ with respect to the basis $g_1(b(\uv, t_1)), \ldots, g_d(b(\uv, t_d))$ of $V_K$. We claim that $\lambda_1, \ldots, \lambda_d\in \cO_w$. If not there exists $i\in \{1, \ldots, d\}$ with $w(\lambda_i) <0$ and $w(\lambda_i)\le w(\lambda_j)$ for all $j\in \{1, \ldots, d\}$. We get
\begin{equation*}
\label{coord}
0\,\,\equiv \,\,\lambda_i^{-1}v \, \, \equiv \,\, \sum_{j=1}^d \, \overline{\lambda_j/\lambda_i}\, b_j \quad \mod \fm_w\otimes_k V
\end{equation*}
which contradicts the assumption that $b_1,\ldots , b_d$ is a basis. The  congruence shows as well that $\overline{\lambda_1},\ldots, \overline{\lambda_d}\in k$ are the coordinates of $v$ with respect to $b_1,\ldots , b_d$. From this the assertion follows immediately.\enddemo

\paragraph{\bf Existence of Shintani cocycles} Assume now $k = \bQ$ and $V=F$. We remark that  $F\!K = F\otimes_{\bQ} K$ as well as $\bR\otimes_{\bQ} K$ are integral domains  since $\bQ$ is algebraically closed in $K$. In fact $F\!K$ is a field. We will denote the quotient field of $\bR\otimes_{\bQ} K$ by $\bK$. The valuation $w$ induces a valuation on $\bK$ with residue field $\bR$. Hence by the same recipe \eqref{order}, $w$ and the parameters $t_0, \ldots, t_d$ induce an ordering on $\bK$. For $i=1, \ldots, d$ the embedding $\sigma_i: F\to \bR$ induce an embedding $F\!K\to \bK$ as well as homomorphisms of $\bR$- resp.\ $\bK$-algebras $F_{\infty} \to \bR$ and $F\!K_{\infty} \colon = F\!K\otimes_K \bK = F\otimes_{\bQ} \bK \to \bK$. By abuse of notation these maps will be denoted by $\sigma_i$ too. We put $F\!K_+ =\{x\in F\!K\mid \sigma_i(x) > 0 \, \forall \, i= 1, \ldots, d\}$ and define $F\!K_{\infty, +}$ similarly. Note that $F\!K_{\infty, +} \cap F = F\!K_+ \cap F = F_+$ and that $F\!K_+$ and $F\!K_{\infty, +}$ are stable under addition and multiplication.  

We define a determinant form $\omega$ on $F$ by $\omega(x_1, \ldots, x_d) = \det(\sigma_i(x_j))$. We fix a $\bQ$-basis $\ux = (x_1, \ldots , x_d)$ of $F$ with $x_1 =1$. For $u_1, \ldots, u_d\in F_{\infty, +}^*$ and $y\in F_{\infty}$ we define 
\begin{eqnarray*}
\label{shincocycle} 
&& z(u_1, \ldots, u_d)(y) \,\, = \,\, c(u_1(b(\ux, t_1)), \ldots, u_d(b(\ux, t_d)))(y)\\
&& = \,\, \sign(\omega(b(u_1 \ux, t_0), \ldots, b(u_d \ux, t_d)))\, 1_{C(b(u_1 \ux, t_0), \ldots, b(u_d \ux, t_d))}(y).
\end{eqnarray*}
The function $y\mapsto z(u_1, \ldots, u_d)(y)$ lies in $\cK$. In fact for $a\in K$ with $w(a) >0$ and $i\in \{1, \ldots, d\}$ we have
\[
\sigma_i(b(\ux, a)) \, = \, \sum_{j=1}^d\, \sigma_i(x_j) a^{j-1} \, = \,  1 + \sum_{j=1}^d\, \sigma_i(x_j) a^{j-1} >0.
\]
Hence $b(u_1 \ux, t_0), \ldots, b(u_d \ux, t_d)\in F\!K_{\infty, +}$ and therefore $z(u_1, \ldots, u_d)(y)  \ne 0$ only if $y\in F_{\infty, +}$. On the other hand by (SC3) there exists finitely many disjoint cones $C_1, \ldots, C_m$ such that $z(u_1, \ldots, u_d) = \pm \sum_{i=1}^m 1_{C_i}$. It follows $C_i \subseteq F_{\infty, +}$ for all $i=1, \ldots, m$. Moreover if $u_1, \ldots, u_d\in F_+^*$ then 
$C_1, \ldots, C_m$ are rational cones hence $z(u_1, \ldots, u_d)\in \cK_{\rat}$.

\begin{lemma}
\label{lemma:shin3}
The map $z: (F_{\infty, +}^*)^d \to \cK, (u_1, \ldots, u_d) \mapsto z(u_1, \ldots, u_d)$ is a homogeneous $(d-1)$--cocyle of $F_+^*$.
\end{lemma}

{\em Proof.} Let $u_0, \ldots, u_d\in F_{\infty, +}^*$. Since $b(u_0\ux, t_0), \ldots, b(u_d \ux, t_d)$ are in general position there exists unique non-zero scalars $\lambda_0, \ldots, \lambda_{d-1}\in \bK^*$ with $\sum_{i=0}^d \, \lambda_i  b(u_i\ux, t_i) = -b(u_d \ux, t_d)$. If all $\lambda_i$ were positive it would imply $-b(u_d \ux, t_d) \in F\!K_{\infty, +}$ contradicting $b(u_d \ux, t_d)\in F\!K_{\infty, +}$. It follows $-b(u_d \ux, t_d)\not\in C(b(u_0\ux, t_0), \ldots, b(u_{d-1}\ux, t_{d-1}))$ and therefore by (SC1) and (SC2) that $z$ is a homogeneous $(d-1)$--cocyle. \enddemo

\begin{prop}
\label{prop:shincoc}
Let  $[z_{\rat}]$ be the cohomology class of the cocycle $z_{\rat}: (F_+^*)^d \to \cK_{\rat}, (u_1, \ldots, u_d) \mapsto z(u_1, \ldots, u_d)$. Then either $[z_{\rat}]$ or $-[z_{\rat}]$ is a Shintani cocycle.
\end{prop}

{\em Proof.} Let $g\in \Maps(F_{\infty,+}, \bZ)$ be the image of $[z_{\rat}]$ under \eqref{shinmap2}. Since the image of $[z_{\rat}]$ under the canonical map $H^{d-1}(F_+^*, \cK_{\rat}) \to H^{d-1}(F_+^*, \cK)$ is the restriction (to $F_+^*$) of the cohomology class $[z]\in  H^{d-1}(F_{\infty,+}^*, \cK)$ of $z: (F_{\infty, +}^*)^d \to \cK$ we get $g = \psi_{E_+}([z])$. Therefore by Lemma \ref{lemma:rescor} $g$ is $\sqrt[\infty]{E_+}$-invariant i.e.\ we have $g(y) = g(\ep y)$ for all $\ep\in \sqrt[\infty]{E_+}$ and $y\in F_{\infty, +}$. By (\cite{colmez}, Lemme 2.1 and 2.2) there exists elements $\ep_1, \ldots, \ep_{d-1}$ of $E_+$ with
\begin{itemize}
\item[(i)] The subgroup $E'$ of $E_+$ generated by $\ep_1, \ldots, \ep_{d-1}$ is free of rank $d-1$. 
\item[(ii)] For all $\tau\in S_{d-1}$ put $f_{1,\tau} = 1$ and  $f_{i,\tau}= \prod_{j<i} \ep_{\tau(j)}$ for $2\le i\le d$. Then $(f_{1,\tau}, \ldots, f_{d,\tau})$ is linearly independent and for all $\tau\in S_{d-1}$ we have 
$\sign(\omega(f_{1,\tau}, \ldots, f_{d,\tau})) = \sign(\tau)$.
\item[(iii)] For $\tau\in S_{d-1}$ put $C_{\tau} = C(f_{1,\tau}, \ldots, f_{d,\tau})$. Then $\ep C_{\tau} \cap \ep' C_{\tau'} = \emptyset$ for $(\ep, \tau), (\ep', \tau')\in E'\times S_{d-1}$ with $(\ep, \tau)\ne (\ep', \tau')$.
\end{itemize}
For $\tau\in S_{d-1}$ we also set $\widetilde{C}_{\tau} = C(b(f_{1,\tau} \ux, t_0), \ldots, b(f_{d,\tau} \ux, t_d))\cap F_{\infty, +}$. Since the closure $\overline{C}_{\tau}$ of $C_{\tau}$ in $F_{\infty, +}$ is $C[f_{1,\tau}, \ldots, f_{d,\tau}]$ we have $C_{\tau} \subseteq  \widetilde{C}_{\tau} \subseteq \overline{C}_{\tau}$ by Lemma \ref{lemma:shin2}. By Lemma \ref{lemma:rescor} (a) we get $g = \psi_{E'}([z])$ hence by Remark \ref{remarks:eta} and Lemma \ref{lemma:shin2} we obtain 
\begin{equation*}
g(y)\! = \!\pm \!\sum_{\ep\in E'} \!\left( \sum_{\tau\in S_{d-1}}  \sign(\tau) \,z(f_{1,\tau}, \ldots, f_{d,\tau})\right)\!\!(\ep^{-1}y) \!=\!\pm \!\!\!\!\!\!\! \sum_{\ep\in E', \tau\in S_{d-1}} \!\! 1_{\ep \widetilde{C}_{\tau}}(y).
\end{equation*}
For pairs $(\ep, \tau)\ne (\ep', \tau')$ condition (iii) implies $\ep C_{\tau} \cap \ep' \overline{C}_{\tau'} = \emptyset$ (because $\ep C_{\tau}$ is open). In particular $C \colon = C_{\id} = C(1, f_1, \ldots, f_{d-1})$ does not intersect $\ep \widetilde{C}_{\tau}$ for all pairs $(\ep, \tau)\ne (1, \id)$ and consequently $g$ is constant $ =\pm 1$ on $C$. Hence by Lemma \ref{lemma:rescor} the function $g$ is constant $=\pm 1$ on the set 
\begin{equation*}
\label{offengitter}
\sqrt[\infty]{E_+} \cdot C \,\,=\,\, \{y \in F_{\infty, +}\mid\, y = \ep y_0 \,\, \mbox{for some $\ep\in \sqrt[\infty]{E_+}$ and $y_0\in C$}\,\}.
\end{equation*}
To finish the proof we have to show 
$\sqrt[\infty]{E_+} \cdot C= F_{\infty, +}$ or that the image of $\sqrt[\infty]{E_+} \cdot C$ under $\Log$ is $= \bR^d$. However $\Log(\sqrt[\infty]{E_+})$ is a $\bQ$-vector space which spans $\bR^d_0$ and for any $t\in \bR$ the intersection of the open set $\Log(C)$ with $\bR^d_t\colon = \{z= (z_1,\ldots, z_d)\in \bR^d\mid \,  \sum_{i=1}^d z_i = t\}$ in nonempty. Hence $\bR^d_t\subseteq \Log(\sqrt[\infty]{E_+} \cdot C)$ for all $t\in \bR$ and therefore $\bR^d = \Log(\sqrt[\infty]{E_+} \cdot C)$.
\enddemo

\paragraph{\bf $(S,T)$-Shintani cocysles.} For a finite set $S$ of nonarchimedian places of $F$ we denote by $\cK^S$ (resp.\ $\cK_S$) the 
subgroup of $\cK_{\rat}$ generated by the characteristic functions $1_C$ of Shintani cones $C$ generated by elements in $E^S_+$ (resp.\ $E_{S, +}$).
Note that $\cK^S$ resp.\ $\cK_S$ is a $E^S_+$ resp.\ $E_{S,+}$-stable subspace of $\cK_{\rat}$. For a prime number $q$ we write $\cK_q$ for $\cK_{S_q}$. We have $\cK_S = \underset{\lra}{\lim}_T \cK^T$ where $T$ runs through all finite subsets of $\bfP_F^{\infty}-S$.

Let $S$, $T$ be finite disjoint subsets $\bfP_F^{\infty}$. Consider the composite
\begin{equation}
\label{edgec2sq}
H^{d-1}(E^S_+, \cK_T)\stackrel{\res}{\lra} H^{d-1}(E_+, \cK_T)\stackrel{\cap \eta_{E_+}}{\lra} H_0(E_+, \cK_T).
\end{equation}

\begin{definition}
\label{definition:shintanisq}
(a) A $T$-integral Shintani decomposition $\cA$ is a Shintani set which can be written as a finite disjoint union of Shintani cones each generated by elements in $E_{T,+}$. 

\noi (b) An $(S,T)$-Shintani cocycle is a cohomology class $\alpha\in  H^{d-1}(E^S_+, \cK_T)$ which is mapped under \eqref{shinmap} to the class $[1_{\cA}]\in H_0(E_+, \cK_T)$ of the characteristic function of a {\it $T$-integral Shintani decomposition} $\cA$. If $T= S_q$ for a prime number $q$ then a $(S,S_q)$-Shintani cocycle will be also called $(S,q)$-Shintani cocycle.
\end{definition}

\begin{prop}
\label{prop:shincocsq}
For a finite subset $S$ of $\bfP_F^{\infty}$ there exists another such set $S_0 \supseteq S$ so that an $(S,T)$-Shintani cocycle exists for all $T$ which are disjoint from $S_0$. In particular there exists an $(S,q)$-Shintani cocycle for almost all prime numbers $q$.
\end{prop}

{\em Proof.} Let $\alpha\in  H^{d-1}(F_+^*, \cK_{\rat})$ be a Shintani cocycle and let $\alpha_S$ be its image under $\res: H^{d-1}(F_+^*, \cK_{\rat})\to H^{d-1}(E^S_+, \cK_{\rat})$. Since $E^S_+$ is finitely generated the functor $H^{d-1}(E^S_+, \wcdot)$ commutes with direct limits. Hence there exists a finite set $S_0\supseteq S$ and $\alpha_0\in H^{d-1}(E^S_+, \cK^{S_0})$ such that $\iota_*(\alpha_0) = \alpha_S$ where $\iota:\cK^{S_0}\subseteq \cK_{\rat}$ is the inclusion.  

There exists a Shintani decomposition $\cA$ such that the image $\alpha$ under \eqref{edgec2} is equal to $[1_{\cA}]\in H_0(E_+, \cK_{\rat})$. By enlarging $S_0$ if necessary we may assume that $\cA$ can be written as a finite disjoint union of Shintani cones generated by elements in $E^{S_0}_+$. This may not necessarily imply that the image of $\alpha_0$ under 
\begin{equation*}
\label{edgec2ss0}
H^{d-1}(E^S_+, \cK^{S_0})\stackrel{\res}{\lra} H^{d-1}(E_+, \cK^{S_0})\stackrel{\cap \eta_{E_+}}{\lra} H_0(E_+, \cK^{S_0}) 
\end{equation*}
is equal to the class of $1_{\cA}$. However by further enlarging $S_0$ we may assume this as well (here we use that $H_0(E_+, \wcdot)$ commutes with direct limits). It is now obvious that for all $T$ disjoint from $S_0$ the image of $\alpha_0$ under the canonical map $H^{d-1}(E^S_+, \cK^{S_0})\to H^{d-1}(E^S_+, \cK_T)$ is a Shintani $(S, T)$-cocycle.\enddemo

\section{Integrality properties of $L$-values attached to Shintani cones} 
\label{section:cassoutrick}

\paragraph{\bf Locally constant functions on adeles and ideles} Our aim now is to relate the function space $\cC_c^0(S_1, S_2, R)$ to the Schwartz space $\cS(\bA\!^{\infty}, R)$ i.e.\ the space of compactly supported locally constant functions $A \to R$.

In general for a locally compact totally disconnected topological ring $A$ and a ring $R$ the  Schwartz space $\cS(A,R)$ is defined as $\cS(A, R) = C^0_c(A,R)$.  The group $A^*$ acts on $\cS(A, R)$ by $(af)(x) \colon = f(a^{-1}x)$ for $a\in A^*$, $f\in \cS(A, R)$ and $x\in A$. Using the embedding $\barQ \hookrightarrow \bC_p$ we view $\cS(A, \barQ)$ as a subspace of $\cS(A, \bC_p)$ and denote the induced $p$-adic maximums norm  \eqref{maxnorm} on $\cS(A, \barQ)$ also by $\| \wcdot\|_p$.

In order to relate $\cC_c^0(S_1, S_2, R)$ to $\cS(\bA\!^{\infty}, R)$ we first consider the local 
case. For $v\in \bP^{\infty}_F$, $\phi\in C_c(F_v^*, R)^{U_v}$ and $x\in F_v^*$ the infinite sum 
\begin{equation*}
\label{unitint}
(\sum_{n=0}^{\infty}\, \varpi^n \phi)(x) \colon = \sum_{n=0}^{\infty}\, \phi(\varpi^{-n} x)
\end{equation*}
is finite and one easily checks that $F_v^* \to R, \, x\mapsto (\sum_{n=0}^{\infty}\, \varpi^n \phi)(x)$ extends to a function in $C^0_c(F_v, R)$. For example if $\phi= 1_{U_v}$ then $\sum_{n=0}^{\infty}\, \varpi^n \phi = 1_{\cO_v}$. We obtain a $F_v^*$-equivariant $R$-linear isomorphism 
\begin{equation}
\label{unitint1}
\delta_v: C_c^0(F_v^*,R)^{U_v} \lra \cS(F_v,R)^{U_v}, \,\,\, \phi\mapsto \sum_{n=0}^{\infty}\, \varpi^n \phi
\end{equation}
which is characterized by $\delta_v(1_{xU_v}) = 1_{x\cO_v}$.

Now consider $\cS(\bA\!^{\infty}, R)$ with its canonical $\bI^{\infty}$-action. We have 
\begin{equation*}
\label{unitint2}
\cS(\bA\!^{\infty}, R) \,\,\, \cong \,\,\, {\bigotimes_{v\nmid\infty}}' \, \cS(F_v,R)
\end{equation*}
where the restricted tensor product $\bigotimes'$ is taken with respect to the family of functions $\{1_{\cO_v}\}_v$. Thus by taking the tensor product of the maps \eqref{unitint1} we obtain a canonical $\bI^{\infty}$-equivariant isomorphism 
\begin{equation*}
\label{unitint2a}
\Delta: C_c^0(\bI^{\infty}/U^{\infty},R) \,\,\, =\,\,\, C_c^0(\bI^{\infty},R)^{U^{\infty}} \,\lra \,\cS(\bA\!^{\infty}, R)^{U^{\infty}}.
\end{equation*}
It can be considered as a linearization of the map sending an idele to the corresponding fractional $\cO_F$-ideal. Indeed, it is characterized by $\Delta(1_{yU^{\infty}})(x) = 1_{\fa}(x)$ for $y\in \bI^{\infty}$ and $x\in F$ where $\fa\in \cI$ denotes the ideal corresponding to $yU^{\infty}$ under the isomorphism $\bI^{\infty}/U^{\infty}\cong \cI$.

More generally, if $S_1, S_2$ are disjoint finite subsets of $\bP^{\infty}_F$ and $S = S_1\cup S_2$ then by taking the tensor product of the canonical inclusion $C_c^0(\bA_{S_1}\times \bA_{S_2}^* ,R)\to \cS(\bA_S, R)$ ("extension by zero") with the tensor product of the maps \eqref{unitint1} for $v\not\in S$ we obtain a $\bI^{\infty}$-equivariant monomorphism
\begin{equation*}
\label{unitint3}
\cC_c^0(S_1, S_2 ,R)  \,\lra \, \cS(\bA\!^{\infty}, R)^{U^{S, \infty}}.
\end{equation*}
It maps a function of the form $\phi= \bigotimes \phi_v$ with $\phi_v\in C^0_c(F_v^*, R)^{U_v}$ for all $v\not\in S$ to $\widetilde{\phi} = \bigotimes \widetilde{\phi}_v$ with $\widetilde{\phi}_v= \phi_v$ if $v\in S$ and $\widetilde{\phi}_v= \delta_v(\phi_v)$ if $v\not\in S$.

Similarly, if $q$ is a prime number with $S_q \cap S = \emptyset$ then by the same procedure we obtain a canonical $\bI^{q, \infty}$-equivariant monomorphism
\begin{equation}
\label{unitint4}
\cC_c^0(S_1, S_2, R)^q \,\lra \, \cS(\bA\!^{q, \infty}, R)^{U^{S,q, \infty}}.
\end{equation}

\paragraph{\bf Solomon-Hu pairing} 

For a finite dimensional $\barQ$-vector space $V$ we denote by $\barQ[\![ V]\!]$ the algebra $\prod_{n\ge 0} \Sym^n V$ and let $\barQ(\!( V)\!)$ be its quotient field. A choice of a basis $(v_1, \ldots v_m)$ of $V$ induces isomorphisms between $\barQ[\![ V]\!]$ resp.\ $\barQ(\!( V)\!)$ and the power series ring $\barQ[\![ z_1,\ldots, z_m]\!]$ resp.\ the field of Laurent series $\barQ(\!( z_1, \ldots, z_m)\!)$. Elements in $\barQ[\![ V]\!]$ can be written as formal sums $\sum_{n\ge 0} v_n$ with $v_n\in \Sym^n V$. We denote the augmentation map by
\begin{equation*}
\label{aug}
\ev_0: \,\barQ[\![ V]\!]\lra \barQ, \, \psi= \sum_{n\ge 0} v_n\, \mapsto \, \ev_0(\psi) = v_0
\end{equation*}
(thus if we think of $\psi$ as a power series in $z_1,\ldots, z_m$ then $\ev_0(\psi)$ is the evaluation at $z_1 =\ldots = z_m =0$). We have $\barQ[\![ V]\!]^* = \{\psi\in \barQ[\![ V]\!]\mid \ev_0(\psi)\ne 0\}$.

For $V = F\otimes_{\bQ} \barQ$ we put $\cR = \barQ[\![ F\otimes_{\bQ} \barQ]\!]$ and $\cQ = \barQ(\!( F\otimes_{\bQ} \barQ)\!)$. The field $F$ can be viewed as a subset of $\cR$ via the embedding $\iota: F\hookrightarrow F\otimes \barQ = \Sym^1(F\otimes \barQ)\subset \cR$, i.e.\ $\iota(x) = \sum_{n\ge 0} v_n$ with $v_1 = x$ and $v_n =0$ for all $n\ne 1$. The multiplication in $F$ induces a $F^*$-action on $\cQ$, $F^*\times \cQ\to \cQ, (x, g) \mapsto x\star g$ which is characterized by $x\star (\iota(x_1)\cdots \iota(x_n)) = \iota(xx_1)\cdots \iota(xx_n)$ for $x,x_1,\ldots,x_n\in F$.

Solomon and Hu \cite{solomon-hu} (see also \cite{hill}) have constructed a pairing
\begin{equation}
\label{solhu1}
\llangle\,\,\,, \,\,\, \rrangle\,:\,\, \cK_{\rat} \times \cS(\bA\!^{\infty}, \barQ) \,\lra \,\cQ
\end{equation}
with the following properties

\noi (SH1) $\llangle xf, x\Phi \rrangle = x \star \llangle f, \Phi \rrangle$ for all $x\in F^*$, $f\in \cK_{\rat}$ and $\Phi\in \cS(\bA\!^{\infty}, \barQ)$.

\noi (SH2) If $C$ is a Shintani cone $C = C(x_1,  \ldots x_m)$ with $x_1, \ldots x_m\in F_+$ linearly independent over $\bQ$ and if $\Phi\in \cS(\bA\!^{\infty}, \barQ)$ is invariant under translation by $\sum_{i=1}^m \bZ x_i$ then \
\begin{equation*}
\label{solhu3}
\llangle 1_C, \Phi \rrangle\ \, = \, \prod_{i=1}^m (1-\exp(x_i))^{-1}\,\cdot\, \sum_{x\in F\cap P(x_1, \ldots, x_m)} \Phi(x) \exp(x)
\end{equation*}
where $\exp(y) = \sum_{i=0}^{\infty} \iota(y)^n/n!\in \cR$ and
\begin{equation*}
P(x_1, \ldots, x_m) \,\, =\,\, \{\sum_{i=1}^m \, t_i x_i\mid \, t_i\in \bR, 0 < t_i \le 1 \,\,\forall \, i=1, \ldots ,m\}.
\end{equation*}

\paragraph{\bf Cassou-Nogu{\`e}s' trick} For a Shintani cone $C$ and $\Phi\in \cS(\bA\!^{\infty}, \barQ)$ we consider the Dirichlet series 
\begin{equation*}
\label{Dirichlet}
L(\Phi, C; s) \,=\, \sum_{x\in F} \, \Phi(x)  1_C(x) \Norm(x)^{-s} .
\end{equation*}
Following Cassou-Nogu{\`e}s \cite{cassou} we study its value at $s=0$ for 
 certain $C$ and $\Phi$.
 
Let $q$ be a prime number $\ne p$ and let $\fq$ be a place of $F$ above $q$. We define $\phi_{\fq}\in \cS(F_{\fq}, \barQ)$ by $\phi_{\fq} = 1_{\cO_{\fq}} -  \Norm(\fq) \,(\varpi_{\fq}1_{\cO_{\fq}})$ i.e.\
\begin{equation*}
\phi_{\fq}(x)  \,\,\, = \,\,\, \left\{\begin{array}{ll} 1 &  \mbox{if $x\in U_{\fq}$,}\\
1-\Norm(\fq) &  \mbox{if $x\in \fq\cO_{\fq}$,}\\
0 &  \mbox{otherwise.}\\
\end{array}\right.
\end{equation*}
Let $\psi_{\fq}: F_{\fq}\to \barQ^*$ denote an (additive) character with $\Ker(\psi_{\fq}) = \cO_{\fq}$ and define $\langle x,y\rangle_{\fq} \colon = \psi_{\fq}(xy)$ for $x,y\in F_{\fq}$. One can easily see that $\phi_{\fq}$ is the Fourier transform of the function $1_{\varpi_{\fq}^{-1} U_{\fq}}\in \cS(F_{\fq}, \barQ)$. Hence by the Fourier inversion formula we obtain
\begin{equation}
\label{fourier}
\phi_{\fq}(x) \,\,\, =\,\,\, -\sum_{y\in Q} \, \langle x,y\rangle_{\fq}
\end{equation}  
for all $x\in F_{\fq}$ where $Q$ is a system of representatives of the set of cosets $\varpi_{\fq}^{-1} \cO_{\fq}/\cO_{\fq} -\{\cO_{\fq}\} = \{y + \cO_{\fp}\mid y\in \fq^{-1} \cO_{\fq}-\cO_{\fq}\}$. Note that if $\Norm(\fq) = q$ then $\langle x,y\rangle_{\fq}$ is a primitive $q$-th root of unity for all $y\in \cQ$ and $x\in U_{\fq}$.  We define $\phi_q\in \cS(F_q, \barQ)$ as the (tensor) product of $\phi_\fq$ and the functions $1_{\cO_v}$ for $v\mid q, v\ne \fq$
\begin{equation}
\label{specialcassou}
\phi_q \,\, =  \,\, \phi_{\fq} \otimes \bigotimes_{v\mid q, v\ne \fq} \, 1_{\cO_v}\in \bigotimes_{v\mid q} \cS(F_v, \barQ) \,\, \cong \,\, \cS(F_q, \barQ).
\end{equation}

Assume now that $\Norm(\fq) = q$ (for example we may choose a prime $q\ne p$ which splits completely in $F$ and take any place $\fq$ above $q$). In the following Lemma we identify $\cS(\bA\!^{q, \infty}, \barQ)\otimes \cS(F_q, \barQ)$ with $\cS(\bA\!^{\infty}, \barQ)$. Thus for $\phi \in \cS(\bA\!^{q, \infty}, \barQ)$ we regard $\phi\otimes \phi_q$ as an element of $\cS(\bA\!^{\infty}, \barQ)$ (the function $\phi_q$ is defined by \eqref{specialcassou}).

\begin{lemma}
\label{lemma:lat0}
Let $\phi \in \cS(\bA\!^{q, \infty}, \barQ)$ and let $C = C(x_1,  \ldots, x_m)$ be a Shintani cone generated by $x_1, \ldots, x_m \in E_{q,+}$. 

\noi (a) We have $\llangle 1_C, \phi\otimes \phi_q \rrangle\in \cR$. Moreover,
\begin{equation}
\label{pmeas}
|\ev_0(\llangle 1_C, \phi\otimes \phi_q \rrangle)|_p \,\,\,  \le \,\,\, \|\phi\|_p.
\end{equation}

\noi (b) $L(\phi\otimes \phi_q, C; s)$ converges absolutely for $\Real(s) > m/d$ and extends holomorphically to the whole complex plane. At $s=0$ we have 
\begin{equation*}
\label{Lvalue}
L(\phi\otimes \phi_q, C;0)\,\, =\,\, \ev_0(\llangle 1_C, \phi\otimes \phi_q \rrangle).
\end{equation*}
\end{lemma}

{\em Proof.} Put $\Phi = \phi\otimes \phi_q$ and $\Phi_0 = \phi\otimes \bigotimes_{v\mid q} 1_{\cO_v}\in \cS(\bA\!^{\infty}, \barQ)$. By \eqref{fourier} we get
\begin{equation}
\label{fourier2}
\Phi(x) \,\,\, =\,\,\, -\sum_{y\in Q} \, \langle x,y\rangle_{\fq} \Phi_0(x)
\end{equation}  
for all $x\in F$ (we consider elements of $\cS(\bA\!^{\infty}, \barQ)$ also as functions on $F$ via the diagonal embedding $F\hookrightarrow \bA\!^{\infty}$). 
There exists two fractional $\cO_F$-ideals $\fa\subseteq \fb\subseteq F$ such that $\Phi_0$ has support in $\fb$ and is constant modulo $\fa$, i.e.\ $\supp(\Phi_0)\cap F \subseteq \fb$ and $\Phi_0(x+ a) = \phi(x)$ for all $a\in \fa$ and $x\in \bA\!^{\infty}$. Since the $q$-component of $\Phi_0$ is $=  \bigotimes_{v\mid q} 1_{\cO_v}$ we may assume that no prime of $\cO_F$ above $q$ occurs in the prime decomposition of $\fa$ and $\fb$. 

Fix $i\in \{1,\ldots, m\}$. Since $x_i$ lies in $E_q$ there exists a positive integer $M$ prime to $q$ such that $M x_i\in \fa$. 
Indeed, no prime ideal above $q$ occurs in the prime decomposition of $x_i^{-1} \fa$ and we can choose $M$ to be the positive generator of $x_i^{-1} \fa\cap \bZ$.

By replacing each $x_i$ by some multiple $Mx_i$ with $q\nmid M$ we can (and will) assume that $x_1, \ldots , x_m\in \fa$, i.e.\ $\Phi_0$ is invariant under translation by $\sum_{i=1}^m \bZ x_i$. Since $x_i$ is a unit at $\fq$ (i.e.\ $\ord_{\fq}(x_i)=0$) we remark that $\langle x_i,y\rangle_{\fq}$ is a primitive $q$-th root of unity for all $y\in Q$ and  $i=1,\ldots, m$. So by \eqref{fourier2} the function $\Phi$ is invariant under translation by $\sum_{i=1}^m \bZ (q x_i)$. 

Since $\wP= P(qx_1, \ldots, qx_m)$ is the disjoint union of sets of the form $(\sum_{i=1}^m  n_i x_i)+P$ with $P= P(x_1, \ldots, x_m)$ and $0\le n_1,\ldots, n_m\le q-1$ and because $\Phi_0$ is constant modulo $\sum_{i=1}^m  \bZ x_i$ we obtain using \eqref{fourier2}
\begin{eqnarray*}
&& \sum_{x\in F\cap \wP} \Phi(x) \exp(x)\,\,  =\,\, - \sum_{y\in Q}\sum_{x\in \fb\cap \wP}  \langle x,y\rangle_{\fq}\Phi_0(x) \exp(x)\\
&& = \,\,  - \sum_{y\in Q} \sum_{x\in \fb\cap P}  \langle x,y\rangle_{\fq}\Phi_0(x)\exp(x)\sum_{n_1,\ldots, n_m=0}^{q-1} \, \prod_{i=1}^m  (\langle x_i,y\rangle_{\fq}\exp(x_i))^{n_i}\\
&& = \,\,  - \prod_{i=1}^m (1- \exp(qx_i))\wcdot \sum_{y\in Q} \frac{\sum_{x\in \fb\cap P}  \langle x,y\rangle_{\fq}\Phi_0(x)\exp(x)
}{\prod_{i=1}^m(1- \langle x_i,y\rangle_{\fq}\exp(x_i))}.
\end{eqnarray*}
By (SH2) we obtain
\begin{eqnarray*}
&&\llangle 1_C, \Phi \rrangle \,\,\, =\,\,\, \prod_{i=1}^m (1-\exp(q x_i))^{-1}\,\cdot\, \sum_{x\in F\cap \wP} \Phi(x) \exp(x)\\
&& = \,\,  - \sum_{y\in Q} \frac{\sum_{x\in \fb\cap P}  \langle x,y\rangle_{\fq}\Phi_0(x)\exp(x)
}{\prod_{i=1}^m(1- \langle x_i,y\rangle_{\fq}\exp(x_i))}\,\, \in \, \cR
\end{eqnarray*}
since $\ev_0(1- \langle x_i,y\rangle_{\fq}\exp(x_i)) = 1-\langle x_i,y\rangle_{\fq}\ne 0$ hence $1- \langle x_i,y\rangle_{\fq}\exp(x_i)\break \in \cR^*$ for all $y\in Q$ and $i=1,\ldots, m$. Moreover we get
\begin{equation*}
\ev_0(\llangle 1_C, \Phi \rrangle)\,\,\, = \,\,\, - \sum_{y\in Q} \,\frac{\sum_{x\in \fb\cap \cP}  \langle x,y\rangle_{\fq}\Phi_0(x)}{\prod_{i=1}^m\, (1 - \langle x_i,y\rangle_{\fq})}.
\end{equation*}
To deduce \eqref{pmeas} note that $|1 -  \langle x_i,y\rangle_{\fq}|_p =1$ hence 
\[
|\ev_0(\llangle 1_C, \Phi \rrangle)|_p \,\, \le \,\, \max_{y\in \cR} |\sum_{x\in \fb\cap \cP}  \langle x,y\rangle_{\fq}\Phi_0(x)|_p \,\, \le \,\,\| \Phi_0\|_p = \|\phi\|_p.
\]
(b) For $L(\Phi, C; s)$ we obtain
\begin{eqnarray*}
&& L(\Phi, C; s) \,\,\, =\,\,\, - \sum_{y\in Q} \, \sum_{x\in \fb\cap C}\, \langle x,y\rangle_{\fq}\Phi_0(x) \Norm(x)^{-s}\\
&& = \,- \sum_{y\in Q} \sum_{x\in \fb\cap \cP} \langle x,y\rangle_{\fq}\Phi_0(x) 
 \sum_{n_1,\ldots, n_m=0}^{\infty} \prod_{i=1}^m \langle x_i,y\rangle_{\fq}^{n_i}\Norm(x + \sum_{i=1}^m n_i x_i)^{-s}.
\end{eqnarray*}
Fix $y\in Q$ and put $\xi_i = \langle x_i,y\rangle_{\fq}$. To deduce (b) it is enough to show that the Dirichlet series 
\[
 \sum_{n_1,\ldots, n_m=0}^{\infty} \prod_{i=1}^m  \xi_i^{n_i} \Norm(x + \sum_{i=1}^m n_i x_i)^{-s}
\]
extends to a holomorphic function on the whole complex plane and that its value at $s=0$ is equal to $\frac{1}{\prod_{i=1}^m\, (1 -  \xi_i)}$. This is well-known (see \cite{shintani}, Prop.\ 9 or \cite{cassou}, Thm.\ 5 and Thm.\ 13). \enddemo

\paragraph{\bf $p$-adic measures attached to Hecke characters and cones} 

Let $\chi: \bI/F^* \to \bC^*$ be a Hecke character of finite order with conductor $\ff(\chi)$. 
Our aim is to construct a $E_{q, +}$-equivariant pairing
\begin{equation}
\label{solhu4}
\llangle\,\,\,, \,\,\, \rrangle_{\chi,\fq}\,:\,\, \cK_q \times \cC_c^0(S_1, S_2,  \barQ)^q \,\lra \barQ
\end{equation}
(a variant of the Solomon-Hu pairing) with the property that for fixed $f\in \cK_q$ the map $\llangle f,  \wcdot \rrangle_{\chi,\fq}: \cC_c^0(S_1, S_2, \barQ)^q \to \bC_p, \, \phi\mapsto \llangle f, \phi \rrangle_{\chi,\fq}$ is an element of $\cD^b(S_1, S_2, \bC_p)$.

Let $S_0$ be the set of all $v\in \bfP_F^{\infty}$ which divide $p\ff(\chi)$. We decompose $S_0$ into disjoint sets $S_0 = S_1 \overset{\cdot}{\cup}  S_2 \overset{\cdot}{\cup}  S_3$ where $S_1= \{v\in S_p\mid \chi_v =1\}$, $S_2 = S_p- S_1$ and $S_3 = S- S_p$. Note that $\chi^{S_0}: \bI^{S_0, \infty}\to \barQ^*$ factors through $\bI^{S_0, \infty}/U^{S_0, \infty}\cong \cI^{S_0}$. Hence we may view $\chi^{S_0}$ as a character 
\begin{equation*}
\label{chiideals}
\chi = \chi^{S_0}: \,\cI^{S_0} \,\,\lra \,\, \barQ^*, \,\, \fa \mapsto \chi^{S_0}(\fa).
\end{equation*}

We fix a place $\fq$ of $F$ such that $q =\Norm(\fq)$ is a prime number with $S_q\cap S_0 = \emptyset$. Since $\bI^{q,\infty}/E_{q,+} U^{S_0, q, \infty} \cong \bI/F^*U^{S_0}$ we can (and will) also regard $\chi$ as a character of 
\begin{equation}
\label{qchi}
\chi: \,\bI^{q,\infty} \,\,\lra \,\, \barQ^*
\end{equation}
with $E_{q,+}U^{S_0,q, \infty} \subseteq \Ker(\chi) $. 

Note that since $\chi_v = 1$ for all $v\in S_1$ we can extend \eqref{qchi} to a multiplicative map $\chi: \bA_{S_1}\times \bI^{S_1,q, \infty} \to \barQ^*$.
We define a map 
\begin{equation*}
\label{delta1}
\Delta^{\chi}\,=\, \Delta^{\chi}_{\fq}: \,\, \cC^0_c(S_1, S_2, \barQ)^q \,\lra \,\cS(\bA\!^{\infty}, \barQ)
\end{equation*}
as the composition
\begin{eqnarray*}
\label{delta1a}
&& \cC^0_c(S_1, S_2, \barQ)^q \, \stackrel{\incl}{\hookrightarrow} \, \cC_c^0(S_1, S_2\cup S_3, \barQ)^q \stackrel{\chi\cdot}{\lra} \cC_c^0(S_1, S_2\cup S_3, \barQ)^q \\
&& \hspace{2cm}\,\stackrel{\eqref{unitint4}}{\lra}\, \cS(\bA\!^{q, \infty},\barQ) \,\stackrel{\wcdot\otimes \phi_q}{\lra} \, \cS(\bA\!^{\infty},\barQ).\nonumber
\end{eqnarray*}
Thus for an element $\phi\in \cC_c^0(S_1, S_2, q, \barQ)$ of the form $\phi = \bigotimes_{v\nmid q\infty} \phi_v$ with
\begin{equation*} 
\phi_v\,\in \, \left\{\begin{array}{ll} C_c^0(F_v, \barQ) &  \mbox{if $v\in S_1$,}\\
C_c^0(F_v^*, \barQ) &  \mbox{if $v\in S_2$,}\\ 
C_c^0(F_v^*, \barQ)^{U_v} &  \mbox{if $v\not\in S_p\cup S_q$}
\end{array}\right.
\end{equation*}
we have $\Delta^{\chi}(\phi) = \bigotimes_{v\nmid \infty} \widetilde{\phi}_v$
where $\widetilde{\phi}_v\in \cS(F_v, \barQ)$ is given by
\begin{equation*}
\widetilde{\phi}_v  \,\,\, = \,\,\, \left\{\begin{array}{ll} \chi_v\phi_v &  \mbox{if $v\in S_0$,}\\
1_{\cO_{\fq}} -  \Norm(\fq) \,1_{\varpi_{\fq}\cO_{\fq}} &  \mbox{if $v=\fq$,}\\
1_{\cO_v} &  \mbox{if $v|q$, $v\ne \fq$,}\\
\delta_v(\chi_v\phi_v) &  \mbox{if $v\not\in S_0\cup S_q$.}
\end{array}\right.
\end{equation*}
Since $\phi_v = 1_{U_v}$ for almost all $v$ we have $\widetilde{\phi}_v
= 1_{\cO_v}$ for almost all $v$.

We remark that 
\begin{equation*}
\label{chiequiv}
\Delta^{\chi}(y\phi)(x) \,\,\, = \,\,\, \chi(y) (\Delta^{\chi}(\phi))(y^{-1} x)
\end{equation*}
for $x\in \bA\!^{\infty}$, $y\in \bI^{q,\infty}$ and $\phi\in \cC_c^0(S_1, S_2, q, \barQ)$. In particular $\Delta^{\chi}$ is $E_{q,+}$-equivariant. We also note that $\|\Delta^{\chi}(\phi)\|_p = \|\phi \|_p$ where $\| \wcdot\|_p$ denote the norm  \eqref{maxnorm} on $\cC^0(S_1, S_2, \barQ)$ and $\cS(\bA\!^{\infty}, \barQ)$ respectively.

Let $y = (y_v)_v\in U_{S_0} \times \bI^{S_0, \infty}$ be an idele whose components $y_v$ at places $v$ above $q$ are all $=1$ (hence we can view $y$ as an element of $\bI^{q,\infty}$). We need an explicite description of $\Delta^{\chi}(1_{yU^{q, \infty}})$. For that let $\fa\in \cI^{S_0}\subseteq \cI$ be the ideal corresponding to $yU^{\infty}$. Then for $x\in F$ we have
\begin{equation}
\label{deltaideal}
\Delta^{\chi}(1_{yU^{q, \infty}})(x) \,\,=\,\,  \chi(x^{-1}\fa)(1_{\fa\cap E_{S_0}}(x)- \Norm(\fq)1_{\fa\fq\cap E_{S_0}}(x)).
\end{equation}

According to Lemma \ref{lemma:lat0} the image of a pair $(f, \Delta^{\chi}(\phi))$ with $f\in \cK_q$ and $\phi\in \cC_c^0(S_1, S_2,\barQ)^q$ under the map \eqref{solhu1} lies in $\cR$. We define \eqref{solhu4} by taking the composition of 
\[
\begin{CD}
\cK_q \times \cC_c^0(S_1, S_2, \barQ)^q\,\,@>\incl\times \Delta^{\chi} >> \cK_{\rat} \times \cS(\bA\!^{\infty},\barQ)
\end{CD}
\]
with \eqref{solhu1} and $\ev_0: \cR\to \barQ$, i.e.\ we have $\llangle  f, \phi\rrangle_{\chi,\fq} = \ev_0(\llangle  f, \Delta^{\chi}(\phi) \rrangle)$ for all $f\in \cK_q$ and $\phi\in \cC_c^0(S_1, S_2,  \barQ)^q$. By Lemma \ref{lemma:lat0} (a) we obtain a $E_{q, +}$-equivariant homomorphism 
\begin{equation*}
\label{solhu5}
\cK_q \,\lra \,\cD^b(S_1, S_2, \bC_p)^q , \, f\mapsto \llangle  f, \wcdot \rrangle_{\chi,\fq}.
\end{equation*}

\paragraph{\bf Shintani decomposition and special $L$-values} 

Let $y^{(1)}\!, \ldots , y^{(h)} \in U_{S_0}\times \bI^{S_0, \infty}$ be ideles whose components at places above $q$ are all $=1$ (hence $y^{(1)}, \ldots , y^{(h)}$ can be regarded as elements $\bI^{q,\infty}$) and such that $y^{(1)}, \ldots , y^{(h)}$ is a system of representatives of $\bI^{q,\infty}/U^{q,\infty}E_{q,+} \cong \Cl^+(F)$ (the narrow class group of $F$). We also consider a $q$-integral Shintani decompostion $\cA$. Recall that this means that $\cA$ has a decomposition $\cA = \overset{\cdot}{\bigcup}_{j\in J} C_j$ where $\{C_j\mid \, j\in J\}$ is a finite collection of Shintani cones $C_j$ which are generated by elements of $E_{q,+}$. 

\begin{lemma}
\label{lemma:Lvalue}
Let $\cF = \bigcup_{i=1}^h y^{(i)} U^{q, \infty}$. Then, 
\begin{equation*}
\label{Lvalue1}
\llangle 1_{\cA}, 1_{\cF} \rrangle_{\chi^{-1},\fq}\,\,= \,\, (1- \chi(\fq) \Norm(\fq)) L_{S_p}(\chi, 0).
\end{equation*}
\end{lemma}

{\em Proof.} Recall that for $\Real(s) >1$ we have 
\begin{equation*}
\label{lchi1a}
L_{S_p}(\chi, s) \,=\, \prod_{v\nmid p\infty} L_v(\chi, s) \,=\, \prod_{v\nmid \infty, v\not\in S} L_v(\chi, s)\,=\, \sum_{\fb\in \cI^S, \fb\subseteq \cO_F} \chi(\fb) \Norm(\fb)^{-s}
\end{equation*}
(since $L_v(\chi, s) = 1$ for all $v\mid \ff(\chi)$). Let $\fa_1, \ldots,\fa_h\in \cI^{S_0}$ be the ideals corresponding to $y^{(1)}U^{\infty}, \ldots , y^{(h)}U^{\infty}$. For $i\in \{1, \ldots, h\}$ put $\Phi_i = \Delta^{\chi^{-1}}_{\fq}\!(1_{y^{(i)} U^{q, \infty}})$ and 
\begin{equation*}
\label{lchi2}
L(\Phi_i, 1_{\cA}; s) \,\, =\,\,  \sum_{x\in F} \, \Phi_i(x)  1_{\cA}(x) \Norm(x)^{-s} \,\, = \,\, \sum_{j\in J} L(\Phi_i, C_j,0).
\end{equation*}
So by Lemma \ref{lemma:lat0} the function $L(\Phi_i, 1_{\cA}; s)$ is entire and by \eqref{deltaideal} we have 
\begin{eqnarray*}
\label{lchi3c}
L(\Phi_i, 1_{\cA}; s) & = & \sum_{x\in \fa_i\cap E_{S_0} \cap \cA} \, \chi(x\fa_i^{-1})\Norm(x)^{-s} \\
& & - \chi(\fq)\Norm(\fq) \sum_{x\in \fq\fa_i\cap E_{S_0} \cap \cA} \, \chi(x(\fq\fa_i)^{-1})\Norm(x)^{-s}.
\end{eqnarray*}
Since $x \mapsto x\fa_i^{-1}$ resp.\ $x\mapsto x(\fq\fa_i)^{-1}$ induces bijections 
\begin{eqnarray*}
\fa_i\cap E_{S_0} \cap \cA & \lra & \{\fb\in \cI^{S_0}\mid \fb\sim \fa_i^{-1}, \, \fb\subseteq \cO_F\},\\ \fq\fa_i\cap E_{S_0} \cap \cA & \lra & \{\fb\in \cI^{S_0}\mid \fb\sim (\fq\fa_i)^{-1}, \, \fb\subseteq \cO_F\}
\end{eqnarray*}
we get
\begin{eqnarray*}
\label{lchi3a}
L(\Phi_i, 1_{\cA}; s) & = & \Norm(\fa_i)^s \, (\sum_{\fb\in \cI^S, \fb\subseteq \cO_F, \fb\sim \fa_i^{-1}} \, \chi(\fb) \Norm(\fb)^{-s})\\
&&   - \chi(\fq)\Norm(\fq)^{1+s} \Norm(\fa_i)^s \, (\sum_{\fb\in \cI^S, \fb\subseteq \cO_F, \fb\sim (\fq\fa_i)^{-1}} \, \chi(\fb) \Norm(\fb)^{-s}).
\end{eqnarray*}
In particular for $s=0$ we obtain 
\begin{equation*}
\label{lchi1b}
\sum_{i=1}^h  \, L(\Phi_i, 1_{\cA}; 0) \,=\, (1 -  \chi(\fq)\Norm(\fq)) L_{S_p}(\chi, 0).
\end{equation*}
On the other hand by Lemma \ref{lemma:lat0} $\ev_0(\llangle  1_{C_j}, \Phi_i  \rrangle) = L(\Phi_i, C_j, 0)$ for all $j\in J$ hence
\[
\llangle 1_{\cA}, 1_{\cF} \rrangle_{\chi^{-1},\fq} \,=\, \sum_{i=1}^h \sum_{j\in J} \ev_0(\llangle   1_{C_j} ,\Phi_i \rrangle)\,=\, \sum_{i=1}^h  \, L(\Phi_i, 1_{\cA}; 0)
\]
so the assertion follows. \enddemo

\section{$L_p(\chi, s)$ via cohomology and proof of the main result}
\label{section:padiclseries}

\paragraph{\bf Interpolation property}

As in last section $\chi: \bI/F^* \to \bC^*$ denotes a Hecke character of finite order, $S_1$ the set of places $v\in S_p$ with $\chi_v=1$, $S_2 = S_p-S_1$ and $S_3$ the set places which do not lie above $p$ and divide $\ff(\chi)$.
We choose $y^{(1)}, \ldots , y^{(h)} \in U_{S_0}\times \bI^{S_0, \infty}$ such that  $y^{(1)}, \ldots , y^{(h)}$ is a system of representatives of $\bI^{\infty}/U^{\infty}F_+^* \cong \Cl^+(F)$ (as before $S_0 = S_1 \overset{\cdot}{\cup}  S_2 \overset{\cdot}{\cup} S_3$). Then there exists a finite subset $S_4$ of $\bfP_F^{\infty}$ which is disjoint from $S_0$ and such that $y^{(1)}, \ldots , y^{(h)}$ already lie in $\bI_{S_4}\times U^{S_4, \infty}$. Put $S = S_0 \cup S_4$ so that $\bI^{S, \infty} = U^{S, \infty}F_+^*$.

By Prop.\ \ref{prop:shincocsq} there exists a prime number $q$ such that (i) $S\cap S_q= \emptyset$, (ii) $q$ splits completely in $F$ and (iii) there exists a $(S,q)$-Shintani cocycle $\alpha\in  H^{d-1}(E^S_+, \cK_q)$. We choose a place $\fq$ of $F$ above $q$ and put $T = S- S_p = S_3 \cup S_4$. In the following we embed $\cC^0_c(S_p, \barQ)_T$ into $\cC^0_c(S_p, \barQ)^q$ and the latter into $\cC^0_c(S_p, \barQ)$ via the map \eqref{shapiro4} and \eqref{shapiro2} respectively. So we have 
\[
\cC^0_c(S_p, \barQ)_T\,\, \subseteq \,\, \cC^0_c(S_p, \barQ)^q\,\, \subseteq \,\, 
\cC^0_c(S_p, \barQ).
\]
When restricting \eqref{solhu4} to $\cK_q \times \cC_c^0(S_1, S_2,  \barQ)_T$ we get a $E^S_+$-equivariant pairing
\begin{equation}
\label{solhu6}
\cK_q \times \cC_c^0(S_1, S_2,  \barQ)_T \,\,\lra \,\,\barQ
\end{equation}
which induces a $E^S_+$-equivariant homomorphism 
\begin{equation}
\label{solhu7}
\cK_q \,\lra \,\cD^b(S_1, S_2, \bC_p)_T .
\end{equation}
Consider the map
\begin{eqnarray}
\label{solhu8}
H^{d-1}(E^S_+, \cK_T) & \lra & H^{d-1}(E^S_+,\cD^b(S_1, S_2, \bC_p)_T)\\
& \stackrel{\eqref{shapiro1}}{\lra} & H^{d-1}(F_+^*,\cD^b(S_1, S_2, \bC_p))\nonumber
\end{eqnarray}
where the first arrow is induced by \eqref{solhu7}. We denote the image of $\alpha$ under  \eqref{solhu8} by $\kappa_{\chi,\fq}$ and let $\mu_{\chi,\fq} = \mu_{\kappa_{\chi,\fq}}$ denote the corresponding $p$-adic measure on $\cG_p$ defined by \eqref{cap2a}.

\begin{prop}
\label{prop:intprop}
For all characters $\eta:  \cG_p \to \bQ^*$ we have 
\begin{eqnarray*}
\int_{\cG_p} \,\eta(\gamma) \mu_{\chi, \fq}(d\gamma) & = & (1- (\chi\eta)^{-1}(\fq) \Norm(\fq)) L_{S_p}((\chi\eta)^{-1}, 0).
\end{eqnarray*}
\end{prop}

{\em Proof.} The pairing \eqref{solhu6} when restricted to the subgroup $\cC_c^0(S_p,  \barQ)_T \break \subseteq \cC_c^0(S_1, S_2,  \barQ)_T$ yields a pairing
\begin{equation*}
\label{solhu9}
\cap_{\solhu} \,\,: \,\,H^{d-1}(E^S_+,\cK_q) \times H_{d-1}(E^S_+,\cC^0_c(S_p,\barQ)_T) \,\lra \,\bC
\end{equation*}
and we have 
\begin{equation}
\label{Delta8}
\kappa_{\chi,\fq} \cap \beta \,\, =\,\, \alpha \cap_{\solhu} \beta\qquad \forall \, \, \beta \in H_{d-1}(E^S_+,\cC^0_c(S_p,\barQ)_T).
\end{equation}
Moreover \eqref{solhu6} induces a pairing 
\begin{equation}
\label{solhu10}
H_0(E_+, \cK_q) \times H^0(E_+,\cC^0_c(S_p, \barQ)_T) \,\lra \,\bC
\end{equation}
and the following diagram commutes
\begin{equation}
\label{Delta10}
\begin{CD}
H^{d-1}(E^S_+,\cK_q) @. \,\,\times\,\, @. H_{d-1}(E^S_+,\cC^0_c(S_p, \barQ)_T) @> \cap_{\solhu} >> \bC\\
@VV\eqref{edgec2sq}V @. @AAA @VV \id V\\
H_0(E_+, \cK_q) @. \,\,\times\,\, @. H^0(E_+,\cC^0_c(S_p,  \barQ)_T) @> \eqref{solhu10} >>\bC
\end{CD}
\end{equation}
where the second vertical map is the composite of 
\begin{equation*}
\label{dualedge1}
 H^0(E_+,\cC^0_c(S_p,  \barQ)_T)\,\,\stackrel{\cap \eta_{E_+}}{\lra} \,\,H_{d-1}(E_+,\cC^0_c(S_p, \barQ)_T) 
 \end{equation*}
with $\cor: H_{d-1}(E_+,\cC^0_c(S_p, \barQ)_T)\to H_{d-1}(E^S_+,\cC^0_c(S_p, \barQ)_T)$. 

Put $\cF = \bigcup_{i=1}^h y^{(i)} U^{q, \infty}$ and let $\eta:  \cG_p \to \barQ^*$ be a character. We view $\eta$ as a Hecke character $\bI^{\fq, \infty}/U^{p, \fq, \infty}\to \barQ^*$ (more precisely we denote the composite
 \[
\bI^{q, \infty}/U^{p, q, \infty}\stackrel{\pr}{\lra}\bI^{q, \infty}/E_{\fq,+}U^{p, q, \infty} \cong \bI/F^*U^p \stackrel{\rho}{\lra} \cG_p \stackrel{\eta}{\lra} \barQ^*
\]
also by $\eta$ where $\rho$ is the reciprocity map). Then $\eta \cdot 1_{\cF}\in \cC^0_c(S_p, \barQ)^q$ actually lies in $H^0(E_+,\cC^0_c(S_p, \barQ)_T)$ 
and its class in $H_0(\Gamma, H^0(E_+, \cC^0_c(S_p, \barQ)))$ is mapped under \eqref{recipr1} to $\eta$ (viewed here as a continuous map $\cG_p\to \barQ$). 

Let $\cA$ be a $q$-integral Shintani decomposition such that the $(S,q)$-Shintani cocycle $\alpha$ is mapped to the class of $1_{\cA}$ in $H_0(E_+, \cK_q)$ under \eqref{edgec2sq}. Using \eqref{Delta8}, the commutativity of \eqref{Delta10} and Lemma \ref{lemma:Lvalue} we obtain
\begin{eqnarray*}
\label{lchi3b}
& \int_{\cG_p} \,\eta(\gamma) \mu_{\chi, \fq}(d\gamma) \,\, = \,\, \kappa_{\chi,\fq} \cap \partial(\eta) \,\, =\,\, \alpha \cap_{\solhu} \partial(\eta) \,\, =\,\, \llangle 1_{\cA}, \eta 1_{\cF} \rrangle_{\chi,\fq}  \\
& =\,\,  \llangle 1_{\cA}, 1_{\cF} \rrangle_{\chi\eta,\fq}  \,\, =\,\, (1- (\chi\eta)^{-1}(\fq) \Norm(\fq)) L_{S_p}((\chi\eta)^{-1}, 0).
\end{eqnarray*}
\enddemo

\paragraph{\bf Proof of Theorem \ref{theorem:vanishing}} Recall (\cite{ribet}, 4.6) that there exists a $\bC_p$-valued $p$-adic measure $\mu$ on $\cG_p$ such that for its $\Gamma$-transform we have
\[
L_p(\mu,s) \,\, =\,\, (1- \chi(\fq)\langle \gamma\rangle^{1-s}) \, L_p(\chi,s)
\]
and such that 
\begin{eqnarray*}
\int_{\cG_p} \,\eta(\gamma) \mu(d\gamma) & = & (1- (\eta\chi)(\fq)\Norm(\fq)) L_{S_p}(\eta\chi, 0)\\
\end{eqnarray*}
for all characters $\eta: \cG_p\to \barQ^*$. Since the latter property determines $\mu$ uniquely we deduce from Prop.\ \ref{prop:intprop} that $\mu =\tau_*(\mu_{\chi, \fq})$ where $\tau: \cG_p \to \cG_p, \gamma \mapsto \gamma^{-1}$. Hence 
\[
L_p(\chi,s)\,\, =\,\,  (1- \chi(\fq)\langle \gamma\rangle^{1-s})^{-1}  \, L_p(\mu_{\chi,\fq}, -s)
\]
in a neighborhood of $s=0$ and therefore 
\[
\ord_{s=0} L_p(\chi,s) \,\, =\,\, \ord_{s=0} L_p(\mu_{\chi,\fq}, s) \,\, \ge \,\, \sharp{S_1}
\]
by Theorem \ref{theorem:abstrezc2}. \enddemo


\begin{thebibliography}{XXXX}

\bibitem{barsky} D.\ Barsky, {\it Fonctions zeta $p$-adiques d'une classe de rayon des corps totalement r{\'e}els}. Groupe d'{\'e}tude analyse ultram{\'e}trique 1977-1978, errata 1978-1979.

\bibitem{burns} D.\ Burns, {\it On derivatives of $p$-adic $L$-series at $s=0$}. Preprint 2012.

\bibitem{cassou} P.\ Cassou-Nogu{\`e}s, {\it Valeurs aux entiers n{\'e}gatifs des fonctions zeta et fonctions zeta $p$-adiques.} Invent.\ Math.\ {\bf 51}, 29--59 (1979).

\bibitem{colmez} P.\ Colmez, {\it R{\'e}sidu en $s=1$ des fonctions zeta $p$-adiques.} Invent.\ Math.\ {\bf 91}, 371--389 (1989).

\bibitem{daschar} P.\ Charollois and S.\ Dasgupta, {\it Integral Eisenstein cocycles on $\GL_n$, I: Szech's cocyle and $p$-adic $L$-functions of Totally Real Fields,} in preparation. 

\bibitem{dasgupta1} S.\ Dasgupta, {\it Shintani Zeta Functions and Gross-Stark Units for Totally Real Fields.} Duke Math.\ J.\ {\bf 143} (2008), 225-279.

\bibitem{dasgupta2} S.\ Dasgupta, {\it Unpublished work,} presented at a conference at the CRL in Barcelona, December 2009.

\bibitem{deligneribet} P.\ Deligne and K.\ Ribet, {\it Values of abelian $L$-functions at negative integers over totally real fields.} Invent.\ Math.\ {\bf 59}, 227--286 (1980).

\bibitem{fedgross} L.J.\ Federer and B.H.\ Gross, {\it Regulators and Iwasawa modules.} Invent.\ Math.\ {\bf 62}, 443--457 (1981).

\bibitem{gross} B.H.\ Gross, {\it $p$-adic $L$-series at $s=0$.} J.\ Fac.\ Sci.\ Univ.\ Tokyo Sect.\ IA Math.\ {\bf 28}, 979--994 (1981).

\bibitem{hill} R.\ Hill, {\it Shintani cocycles on $\GL_n$}, Bull.\ L.M.S. {\bf 39}, 993 - 1004. 

\bibitem{solomon-hu} S.\ Hu and D.\ Solomon, {\it Properties of higher-dimensional Shintani generating functions and cocycles on $\PGL_3(\bQ)$.} Proc.\ L.M.S.\ {\bf 82} (2001)

\bibitem{ribet} K.\ Ribet, {\it Report on $p$-adic $L$-functions over totally real fields.} Ast{\'e}risque {\bf 61}, 177--192 (1979).

\bibitem{shintani} T.\ Shintani, {\it On the evaluation of zeta functions of totally real fields.}  J.\ Fac.\ Sci.\ Univ.\ Tokyo Sect.\ IA Math.\ {\bf 23}, 393--417 (1976).

\bibitem{solomon} D.\ Solomon, {\it Algebraic properties of Shintani's generating 
functions: Dedekind sums and cocycles on $\PGL_2(\bQ)$.} Comp.\ Math.\ {\bf 112}, 333--362 (1998).

\bibitem{ich} M.\ Spie{\ss}, {\it On special zeros of $p$-adic $L$-functions of Hilbert modular forms.} Preprint.

\bibitem{wiles} A.\ Wiles, {\it The Iwasawa conjecture for totally real fields.} Ann.\ of Math.\ {\bf 131} (1990), 493--540.
\end{thebibliography}
\end{document}